\documentclass[smallextended,referee,envcountsect]{svjour3}
\smartqed

\usepackage{amsfonts}
\usepackage{amsmath}
\usepackage{amssymb}
\usepackage{graphicx}
\usepackage{caption}
\newtheorem{lem}{Lemma}
\newtheorem{defin}{Definition}
\makeatletter 
\def\smallunderbrace#1{\mathop{\vtop{\m@th\ialign{##\crcr  
$\hfil\displaystyle{#1}\hfil$\crcr 
\noalign{\kern3\p@\nointerlineskip}%
\tiny\upbracefill\crcr\noalign{\kern3\p@}}}}\limits}

\begin{document}
\title{Solving a fractional parabolic-hyperbolic free boundary problem
which models the growth of tumor with drug application using  finite difference-spectral method}

\author{   Sakine Esmaili\and F. Nasresfahani \and M.R. Eslahchi  }

\institute{
Sakine Esmaili  \\
 sakine.esmaili@modares.ac.ir\\
Farzane Nasresfahani\\
 f.nasresfahani@modares.ac.ir \at
Mohammad Reza Eslahchi, Corresponding author \\
eslahchi@modares.ac.ir\\
Department of Applied Mathematics, Faculty of Mathematical Sciences, Tarbiat Modares University, P.O. Box 14115-134 \\
Tehran, Iran
}
\date{}
\maketitle
\begin{abstract}
\noindent In this paper, a free boundary problem modelling the growth of tumor is considered. The model includes two reaction-diffusion equations modelling the diffusion of nutrient and drug in the tumor and three hyperbolic equations describing the evolution of three types of cells (i.e. proliferative cells, quiescent cells and dead cells) considered in the tumor. Due to the fact that in the  real situation, the subdiffusion of nutrient and drug in the tumor can be found, we have changed the reaction-diffusion equations to the fractional ones to consider other conditions and study a more general and reliable model of tumor growth.  Since it is important to solve a problem to have a clear vision of the dynamic of  tumor growth under the effect of the nutrient and drug, we have solved the fractional  free boundary problem.  We have solved the fractional parabolic equations employing a combination of spectral  and finite difference methods and the  hyperbolic equations are solved using characteristic equation and finite difference method. It is proved that the presented method is unconditionally convergent and stable to be sure that we have a correct vision of  tumor growth dynamic. Finally, by presenting some numerical examples and showing the results,  the theoretical statements are justified.
\end{abstract}
\section{Introduction}
Cancer is one of the most leading causes of death and many different tumor types  in the human body are diagnosed such as Glioblastomas,  phyllodes tumors and so on.
Glioblastomas are the most common
and malignant primary brain tumor, which are very aggressive,  with  the ability to recur
despite extensive treatment \cite{glio1,glio}. Phyllodes tumors are breast tumors and most often are benign (even when they are benign, they can  show recurrence), but in some cases they can be malignant \cite{phyll}. Due to the importance of the treatment of the malignant tumors, scientists including the mathematicians have studied the cancer problem using different tools and techniques. Mathematicians have applied mathematical modelling techniques to model  the various aspects
of cancer dynamics such as avascular and vascular tumor growth, invasion, metastasis and so on.  For instance, the authors of \cite{vas1,vas2} studied the mathematical models of vascular tumors  while the author of  \cite{J.H.}  investigated the mathematical models of  avascular tumors. Owing to the fact that, low concentrations of glucose
and oxygen in the inner regions of spheroids may  contribute to the formation of many types of cell subpopulations such as quiescent-, hypoxic-, anoxic- and necrotic  cells \cite{khai}, therefore some tumor growth models have  divided alive cells  into proliferative  and quiescent cells \cite{J.H.,Y. Tao}.  Since, in vitro results show that in the early stages the solid tumors grow  approximately spherically symmetric  \cite{M. Chen}, in most of the models it is assumed that the tumor grows radially-symmetric. As we mentioned above, treatment of tumors and destroying them is very important, therefore many researchers have investigated the models of tumor growth in which the treatment of tumor is taken into account \cite{ther1,ther2,ther3}. In \cite{ther1},  a mathematical
model for the combined treatment of chemotherapy and radiation in Non-Small Cell Lung Cancer patients is developed to improve the
treatment strategies for future clinical trials. The authors of \cite{ther2} applied a system of nonlinear ordinary differential equations  to  analyse a mathematical model of the treatment of colorectal cancer. In the model  the effects of immunotherapy and chemotherapy on the  tumor cells and cancer stem cells is described. 
  In another study \cite{fansari},  applying  a system of four coupled partial differential equations, the interaction between normal, immune and tumor cells in a tumor
with a chemotherapeutic drug is described.
One of the models with mentioned properties, i.e. in which  three types of cells including proliferative, quiescent and dead cells are considered, the tumor is assumed to grow radially symmetric, and the effect of drug on the treatment of the tumor is also taken into account, is presented in \cite{J.H.}, which is briefly  as follows:
\begin{equation} \label{1****} 
\frac{\partial C}{\partial t}=D_1\frac{1}{r^2}\frac{\partial }{\partial r}\left(r^2\frac{\partial C}{\partial r}\right)-F\left(C,P,Q\right), \ \ \ 0<r<R\left(t\right),\ \ t>0,
\end{equation} 
\begin{equation} \label{GrindEQ__3_2_***} 
\frac{\partial C}{\partial r}\left(r,t\right)=0\ at\ r=0,\ C\left(r,t\right)=\overline{C}(t)\ at\ r=R\left(t\right),\ \ t>0, 
\end{equation} 
\begin{equation} 
 C\left(r,0\right)=C_0(r), ~~ 0\leq r\leq R_0, 
\end{equation} 
\begin{equation} \label{2****} 
\frac{\partial W}{\partial t}=D_2\frac{1}{r^2}\frac{\partial }{\partial r}\left(r^2\frac{\partial W}{\partial r}\right)-G\left(W,P,Q\right),\ \ \  0<r<R\left(t\right),\ \ t>0, 
\end{equation} 
\begin{equation} \label{3****} 
\frac{\partial W}{\partial r}\left(r,t\right)=0\ at\ r=0,\ W\left(r,t\right)=\overline{W}(t)\ at\ r=R\left(t\right),\  \ \ \ t>0, 
\end{equation} 
\begin{equation} 
 W\left(r,0\right)=W_0(r), ~~ 0\leq r\leq R_0, 
\end{equation}
\begin{equation} \label{6****} 
\frac{\partial P}{\partial t}+\dot{v}\frac{\partial P}{\partial r}=g_{11}\left(C,W,P,Q,D\right)P+g_{12}\left(C,W,P,Q,D\right)Q+g_{13}\left(C,W,P,Q,D\right)D, 
\end{equation} 
\begin{equation} \label{GrindEQ__3_6_***} 
\frac{\partial Q}{\partial t}+\dot{v}\frac{\partial Q}{\partial r}=g_{21}\left(C,W,P,Q,D\right)P+g_{22}\left(C,W,P,Q,D\right)Q+g_{23}\left(C,W,P,Q,D\right)D, 
\end{equation} 
\begin{equation} \label{7****} 
\frac{\partial D}{\partial t}+\dot{v}\frac{\partial D}{\partial r}=g_{31}\left(C,W,P,Q,D\right)P+g_{32}\left(C,W,P,Q,D\right)Q+g_{33}\left(C,W,P,Q,D\right)D, 
\end{equation} 
\[0\le r\le R\left(t\right),\ \ t>0,\] 
\begin{equation}\label{Ak6***} 
\frac{1}{r^2}\frac{\partial }{\partial r}\left(r^2 \dot{v}\right)=h\left(C,W,P,Q,D\right),\ \ 0<r\leq R\left(t\right),\ \ t>0, 
\end{equation} 
\begin{equation}
\dot{v}\left(0,t\right)=0,\ \ t>0, 
\end{equation} 
\begin{equation} \label{Ak7***} 
\frac{dR(t)}{dt}=\dot{v}\left(R\left(t\right),t\right),\ \ t>0, 
\end{equation} 
\begin{equation} \label{GrindEQ__3_11_***} 
 P\left(r,0\right)=P_0\left(r\right),~~Q\left(r,0\right)=\ Q_0\left(r\right),\ \ \ D\left(r,0\right)=\ D_0\left(r\right),\ \ \ R(0)=R_0,\ \ \ \ 0\le r\le R_0, 
\end{equation} 
where  $C$ and $W$ are the concentration of  nutrient and drug, respectively.  $P$, $Q$ and $D$ are densities of proliferative cells, quiescent cells and dead cells, respectively and $R(t)$ is the radius of tumor at time $t$.  Also, it's assumed that the initial data satisfies the following conditions
\begin{eqnarray}\label{initialhgjh}
&0\le C_0\left(r\right)\le \overline{C}(0),~ 0\le W_0\left(r\right)\le \overline{W}(0),&~ 0\le r\le R_0,\nonumber\\
&\dfrac{\partial C_0(r)}{\partial r}=0~ at~r=0,~ C_0(R_0)=\overline{C}(0),&\nonumber\\
&\dfrac{\partial W_0(r)}{\partial r}=0~ at~r=0,~ W_0(R_0)=\overline{W}(0),&\\
&P_0\left(r\right)\ge 0,\ \ Q_0\left(r\right)\ge 0,\ \ D_0\left(r\right)\ge 0,&~ 0\le r\le R_0,\nonumber\\
&{P}_0\left(r\right)+Q_0\left(r\right)+D_0\left(r\right)=N,&~ 0\le r\le R_0.\nonumber
\end{eqnarray}
In this model, it is assumed that the nutrient and drug diffuse throughout the tumor with diffusion coefficient $D_1$ and $D_2$, respectively. But in a real situation, subdiffusion of nutrient and drug in the tumor can be found. 
 Therefore, many papers are devoted to studying and solving the  fractional-order mathematical models of  the dynamic of cancers and different methods are employed to deal with fractional models.  For instance,  the application of the homotopy perturbation method 
to two-point boundary-value  problems with fractional-order derivatives of Caputo type is studied in \cite{Ates}.
In  another study \cite{Sohail}, the fractional-order mathematical model involved three Michaelis–
Menten nonlinear terms  and two types of treatments is numerically solved. In \cite{Veeresha},  a fractional-order
cancer chemotherapy effect model in Caputo sense is formulated, which
is studied and analysed by Ansarizadeh et al in \cite{fansari}. q-HATM is
used to solve the system of equations
with a chemotherapeutic drug, describing the interaction
among tumor cells, immune cells, and normal cells in a
tumor. \\
In this paper,  we have considered fractional parabolic equations  to deal with a more reliable model  of tumor growth. 
Then, we have solved the obtained problem, which includes two fractional parabolic equations and three hyperbolic equations, employing a combination of finite difference method and spectral method. Solving the problem enables us to have a clear vision of the dynamic of the tumor under the effect of nutrient and drug. For having an accurate vision, it is of great importance to prove the convergence of the method to be sure that the results are trustable. Therefore, we have also proved the unconditional convergence and stability of the method. Finally, by presenting some numerical examples together with the results, the theoretical statements are justified.

\section{Fractional model of tumor growth}
In this article, it is aimed to solve   a fractional parabolic-hyperbolic free boundary problem modelling the growth of tumor with drug application which consists of two fractional parabolic and three hyperbolic differential equations and one ordinary differential equation which are coupled together. The model with integer derivatives (i.e. the model given in \eqref{1****}--\eqref{initialhgjh}) is presented in \cite{J.H.} but we have changed the problem to a fractional one to consider  the subdiffusion of nutrient and drug.  Moreover, since it is supposed that the tumor grows radially symmetric with free boundary so we can change the domain of model to the following fixed domain
\[\left\{\left(\rho {\rm ,\ } t \right){\rm\mid \ 0\ }\le {\rm \ }\rho {\rm \ }\le {\rm \ 1,\ } t {\rm \ }\ge {\rm \ 0}\right\}.\]     
Using the following change of variables
\[\rho =\frac{r}{R\left(t\right)},~c\left(\rho , t \right){\rm =}C\left(r,t\right),~w\left(\rho , t \right){\rm =}W\left(r ,t \right),~~v(\rho,t)=\dot{v}(r,t),\]  
\begin{equation}\label{Yza}
p\left(\rho , t \right){\rm =}P\left(r,t\right),~q\left(\rho , t \right){\rm =}Q\left(r,t\right),~ d\left(\rho , t \right){\rm =}D\left(r,t\right),
\end{equation} 
and by considering the subdiffusion of nutrient and drug the model becomes
\begin{equation}\label{YAS}
\dfrac{\partial c}{\partial t }=\dfrac{\partial^\alpha }{\partial t^\alpha }\left(\dfrac{D_{1}}{(R(t)\rho)^2}\dfrac{\partial }{\partial\rho }\Big({\rho }^{2}\dfrac{\partial c}{\partial\rho }\Big)\right)+\dfrac{{v}(1,t)\rho}{R(t)}\dfrac{\partial c}{\partial\rho}-f\left(c,p,q\right),~ {\rm \ 0}<\rho<1,\  t {\rm >}0,
\end{equation}
\begin{equation}
\dfrac{\partial c}{\partial\rho }\left(0, t \right){\rm =0, }~c\left({\rm 1,} t \right){\rm =}\overline{c}( t), ~~\  t {\rm >}0,
\end{equation} 
\begin{equation}
c\left(\rho ,0\right){\rm =}c_0\left(\rho \right),~{\rm  \ 0}\le \rho \le {\rm 1,}
\end{equation}
\begin{equation}\label{YAS*}
\dfrac{\partial w}{\partial t }=\dfrac{\partial^\alpha }{\partial t^\alpha }\left(\dfrac{D_{2}}{{(R(t)\rho) }^{2}}\dfrac{\partial }{\partial\rho }\Big({\rho }^{2}\dfrac{\partial w}{\partial\rho }\Big)\right)+\dfrac{{v}(1,t)\rho}{R(t)}\dfrac{\partial w}{\partial\rho}-g\left(w,p,q\right),~ {\rm \ 0}<\rho <1,\  t {\rm >}0,
\end{equation}
\begin{equation}
\dfrac{\partial w}{\partial\rho }\left(0, t \right){\rm =0, }~w\left({\rm 1,} t \right){\rm =}\overline{w}( t), ~~\  t {\rm >}0,
\end{equation}
\begin{equation}\label{YAS1}
w\left(\rho ,0\right){\rm =}w_0\left(\rho \right),~{\rm \ \ 0}\le \rho \le {\rm 1,}
\end{equation}
\begin{eqnarray}\label{pqd}
&\dfrac{\partial p}{\partial t }{\rm +}\smallunderbrace{\dfrac{(v-\rho v(1,t))}{R(t)}}_{\nu(\rho,t)}\dfrac{\partial p}{\partial\rho }{\rm =}
(g_{{\rm 11}}\left(c,w,p,q,d\right)p{\rm +}g_{{\rm 12}}\left(c,w,p,q,d\right)q{\rm +}g_{{\rm 13}}\left(c,w,p,q,d\right)d){\rm , }&\\  
&\dfrac{\partial q}{\partial t }{\rm +}\dfrac{{(v-\rho v(1,t))}}{R(t)}\dfrac{\partial q}{\partial\rho }{\rm =} 
(g_{{\rm 21}}\left(c,w,p,q,d\right)p{\rm +}g_{{\rm 22}}\left(c,w,p,q,d\right)q{\rm +}g_{{\rm 23}}\left(c,w,p,q,d\right)d){\rm , }&\\  
&\dfrac{\partial d}{\partial t }{\rm +}\dfrac{{(v-\rho v(1,t))}}{R(t)}\dfrac{\partial d}{\partial\rho }{\rm =}
(g_{{\rm 31}}\left(c,w,p,q,d\right)p{\rm +}g_{{\rm 32}}\left(c,w,p,q,d\right)q{\rm +}g_{{\rm 33}}\left(c,w,p,q,d\right)d),&\\ 
&{\rm 0}\le \rho \le {\rm 1,\ } t {\rm >}0,& \nonumber
\end{eqnarray}
\begin{equation}\label{YAS2}
~~~~~~~~~~~~~~~~~~~p\left(\rho,0\right)=\ p_0\left(\rho\right),~~~q\left(\rho,0\right)=\ q_0\left(\rho\right),~~~ d\left(\rho,0\right)=d_0\left(\rho\right),~~~~{\rm 0}\le \rho \le {\rm 1,\ }
\end{equation}
\begin{eqnarray} \label{HM345} 
&\dfrac{1}{\rho^2}\dfrac{\partial }{\partial \rho}\left(\rho^2\dfrac{{v}}{R(t)}\right)=h\left(c,w,p,q,d\right),\ \ 0<\rho\leq 1,\ \ t>0,& \nonumber\\
&{v}\left(0,t\right)=0,\ \ t>0,& \\
&\dfrac{dR(t)}{dt}={v}\left(1,t\right),\ \ t>0, \ \ R(0)=R_0,& \nonumber
\end{eqnarray}
where $\dfrac{\partial^\alpha}{\partial\rho}={}_0D^\alpha_t$ is the Riemann-Liouville fractional derivative and $0<\alpha<1$ and \[f\left(c,p,q\right)=K_1\left(c\right)p+K_2\left(c\right)q,~~g\left(w,p,q\right)=K_3\left(w\right)p+K_4\left(w\right)q,\] 
\[g_{11}\left(c,w,p,q,d\right)=\left[K_B\left(c\right)-K_Q\left(c\right)-K_A\left(c\right)-G_1\left(w\right)\right]-\frac{1}{N}\left[K_B\left(c\right)p-K_Rd\right],\] 
\[g_{12}\left(c,w,p,q,d\right)=K_P\left(c\right),~g_{13}\left(c,w,p,q,d\right)=0,\] 
\[g_{22}\left(c,w,p,q,d\right)=-\left[K_P\left(c\right)+K_D\left(c\right)+G_2\left(w\right)\right]-\frac{1}{N}\left[K_B\left(c\right)p-K_Rd\right],\]
\[g_{21}\left(c,w,p,q,d\right)=K_Q\left(c\right),~g_{23}\left(c,w,p,q,d\right)=0,\] \[g_{31}\left(c,w,p,q,d\right)=K_A\left(c\right)+G_1\left(w\right),~~g_{32}\left(c,w,p,q,d\right)=K_D\left(c\right)+G_2\left(w\right),\] 
\[g_{33}\left(c,w,p,q,d\right)=-K_R-\frac{1}{N}\left[K_B\left(c\right)p-K_Rd\right],\] 
\[h\left(c,w,p,q,d\right)=\frac{1}{N}\left[K_B\left(c\right)p-K_Rd\right],\] 
\begin{equation*} 
{\ \ \ \ g}_{ij}\left(c,w,p,q,d\right)\ge 0,\ \ \ \ \ \ i\ne j, 
\end{equation*} 
and the initial data satisfies the following conditions
\begin{eqnarray}\label{initial1}
&0\le c_0\left(\rho\right)\le \overline{c}(0),~ 0\le w_0\left(\rho\right)\le \overline{w}(0),&~ 0\le \rho\le 1,\nonumber\\
&\dfrac{\partial c_0(\rho)}{\partial \rho}=0~ at~r=0,~ c_0(1)=\overline{c}(0),&\nonumber\\
&\dfrac{\partial w_0(\rho)}{\partial \rho}=0~ at~r=0,~w_0(1)=\overline{w}(0),&\\
&p_0\left(\rho\right)\ge 0,\ \ q_0\left(\rho\right)\ge 0,\ \ d_0\left(\rho\right)\ge 0,&~ 0\le \rho\le 1,\nonumber\\
&{p}_0\left(\rho\right)+q_0\left(\rho\right)+d_0\left(\rho\right)=N,&~ 0\le \rho\le 1\nonumber.
\end{eqnarray}
In this model  $G_1(w)$ and $G_2(w)$  are the dead rates of the proliferative cells and quiescent cells due to the drug, respectively.  $K_B(c)$ is the mitosis rate of proliferative cells that is dependent on nutrient level $c$, $K_A(c)$ and $K_D(c)$ are death rates of proliferative cells and quiescent cells, respectively.
  $K_P (c)$ and $K_Q(c)$ are the transferring rate of quiescent cells to proliferative cells and the rate of  transferring proliferative cells to quiescent, respectively. $K_R$ is the constant rate of removing dead cells from the tumor.  In the presented model, $(K_1\left(c\right)p+K_2\left(c\right)q)$ and $(K_3\left(w\right)p+K_4\left(w\right)q)$ show consumption rate  of nutrient and  drug, respectively. $\overline{c}(t)$ and $\overline{w}(t)$ are positive functions  showing   nutrient and  drug supply that the tumor receives from its boundary.
  \noindent In the model (\ref{YAS})--(\ref{HM345}) it is assumed that

\noindent \textbf{A.~}$K_1\left(c\right)$, $K_2\left(c\right)$, $K_3\left(w\right)$, $K_4\left(w\right)$,   $K_A\left(c\right),~K_B\left(c\right),~K_D\left(c\right),~K_P\left(c\right)$, $K_Q(c)$, $G_1(w)$ and $G_2(w)$ are  $C^2$-smooth functions.

\noindent \textbf{B.~}$p_0,{q}_0$ and $d_0$ are non-negative $C^1$-smooth functions on $\left[0,R_0\right].$

\noindent \textbf{C.~}  $c(|x|,\tau)$ and $w(|x|,\tau)$ on $\partial Q^1_T\setminus(B_1\times\{t=T\})$ are non-negative  functions, also $c(|x|,\tau)=\Psi(x,\tau)$ and $w(|x|,\tau)=\Psi_1(x,\tau)$ on $\partial Q^1_T\setminus(B_1\times\{t=T\})$, where $Q^1_T:=\{(x,\tau)\in \mathbb{R}^{3}\times \mathbb{R}:~ |x|< 1,~0<\tau\leq T\}$, $B_{1}=\left\{ x\in {\rm \ \mathbb{R}^3:~}| x|\le 1\right\}{\rm}$  and $\Psi(x,\tau),~\Psi_1(x,\tau)\in C^{2+\alpha,1+\frac{\alpha}{2}}(\overline{ Q^1_T})$ for some $0<\alpha<1$ ($ C^{2+\alpha,1+\frac{\alpha}{2}}(\overline{\Omega}\times[0,T])$ for $0<\alpha<1$ and $\Omega\subset \mathbb{R}^n$ is defined in Appendix).

The main aim of this article is to solve this initial-boundary value problem using the spectral  and finite difference method.

\section{Approximating the solution of the problem}\label{3}
In this section, we  approximate the solution of problem (\ref{YAS})--(\ref{HM345}) for $0\leq\rho\leq1$ and $0\leq t\leq T$. Let $ t_n := n t^*$ ($n = 0, 1,\cdots, M$) be mesh points,  where $ t^* := \frac{T}{M}$ is the time step  and M is
a positive integer. The problem is solved employing spectral method and  the following discretization formula \cite{impfrac} for approximating the time fractional derivative 
\[\dfrac{\partial^\alpha u}{\partial  t^\alpha}(\rho, t_{n})=u(\rho,0)\dfrac{t_{n}^{-\alpha}}{\Gamma(1-\alpha)}+\dfrac{1}{\Gamma(1-\alpha)}\sum_{k=0}^{n-1}\int_{ t_k}^{ t_{k+1}}\dfrac{\partial u(\rho,s)}{\partial s}\dfrac{1}{( t_{n}-s)^\alpha}ds=\]
\[u(\rho,0)\dfrac{t_{n}^{-\alpha}}{\Gamma(1-\alpha)}+\dfrac{1}{\Gamma(1-\alpha)}\sum_{k=0}^{n-1}\dfrac{ u(\rho, t_{k+2})-u(\rho, t_{k+1})}{ t^*}\int_{ t_k}^{ t_{k+1}}\dfrac{1}{( t_{n}-s)^\alpha}ds+E^u_t=\]
\[u(\rho,0)\dfrac{t_{n}^{-\alpha}}{\Gamma(1-\alpha)}+\dfrac{( t^*)^{-\alpha}}{\Gamma(2-\alpha)}\sum_{k=0}^{n-1} \Big(u(\rho, t_{k+2})-u(\rho, t_{k+1})\Big)\Big((n-k)^{1-\alpha}-(n-1-k)^{1-\alpha}\Big)+E^u_t=\]
\[u(\rho,0)\dfrac{t_{n}^{-\alpha}}{\Gamma(1-\alpha)}+\dfrac{( t^*)^{-\alpha}}{\Gamma(2-\alpha)}\sum_{k=0}^{n-1} \Big(u(\rho, t_{n+1-k})-u(\rho, t_{n-k})\Big)\Big((k+1)^{1-\alpha}-(k)^{1-\alpha}\Big)+E^u_t=\]
\begin{equation}\label{frac1}u(\rho,0)\dfrac{t_{n}^{-\alpha}}{\Gamma(1-\alpha)}+\sum_{k=0}^{n-1} a_k\Big(u(\rho, t_{n+1-k})-u(\rho, t_{n-k})\Big)+E^u_ t,
\end{equation}
 where 
  \begin{equation}\label{ak0}
  a_k=\dfrac{(k+1)^{1-\alpha}-(k)^{1-\alpha}}{( t^*)^{\alpha}\Gamma(2-\alpha)}, ~~a_{k+1}\leq a_k,~~\|E_t^u\|_{\infty}<C_1(t^*)^{2-\alpha}.
  \end{equation}
In the following without loss of generality we can suppose that $\overline{c}=\overline{w}=0$ and $c_0(\rho)=w_0(\rho)=0$.  We have also assumed that
\[A_n(\rho):=A(\rho, t_n),~~B_n:=B( t_n).\]From \eqref{frac1} and
(\ref{YAS})--(\ref{YAS1}), we get
\begin{eqnarray}\label{YHM4kk}
  & c_{n+1}-c_{n}+\dfrac{c_{n-1}-c_n}{3}-\dfrac{a'_0D_{1}}{(R_{n+1})^2}\dfrac{1}{{\rho }^{2}}\dfrac{\partial }{\partial\rho }\left({\rho }^{2}\dfrac{\partial c_{n+1}}{\partial\rho }\right)-t^*\dfrac{2(2{v}(1,t_{n})-{v}(1,t_{n-1}))\rho}{3R_{n+1}}\dfrac{\partial c_{n+1}}{\partial\rho}=&\nonumber\\
  &-\displaystyle\sum_{k=0}^{n-1} (a'_k-a'_{k+1})\Big(\dfrac{D_{1}}{(R_{n-k})^2}\dfrac{1}{{\rho }^{2}}\dfrac{\partial }{\partial\rho }\left({\rho }^{2}\dfrac{\partial c_{n-k}}{\partial\rho }\right)\Big)+&\nonumber\\
  &\dfrac{2}{3}\Big(-2t^*f\left(c_{n},p_{n},q_{n}\right)+t^*f\left(c_{n-1},p_{n-1},q_{n-1}\right)\Big)-E^c_ t,&\\
  &0< \rho<1,\ \ 0< t\leq T,&\nonumber\\
&\dfrac{\partial c_{n+1}}{\partial\rho }\left(0 \right){\rm =0, }~c_{n+1}\left({\rm 1}  \right){\rm =}0,~~~~~~\ 0< t\leq T,&\nonumber\\
&c_{0}\left(\rho\right){\rm =}0,~~~~~0\leq \rho\leq1,&\nonumber
\end{eqnarray}
\begin{eqnarray}\label{YHM5kk}
 & w_{n+1}-w_{n}+\dfrac{w_{n-1}-w_n}{3}-\dfrac{a'_0D_{2}}{(R_{n+1})^2}\dfrac{1}{{\rho }^{2}}\dfrac{\partial }{\partial\rho }\left({\rho }^{2}\dfrac{\partial w_{n+1}}{\partial\rho }\right)-t^*\dfrac{2(2{v}(1,t_{n})-{v}(1,t_{n-1}))\rho}{3R_{n+1}}\dfrac{\partial w_{n+1}}{\partial\rho}=&\nonumber\\
  &-\displaystyle\sum_{k=0}^{n-1} (a'_k-a'_{k+1})\Big(\dfrac{D_{2}}{(R_{n-k})^2}\dfrac{1}{{\rho }^{2}}\dfrac{\partial }{\partial\rho }\left({\rho }^{2}\dfrac{\partial w_{n-k}}{\partial\rho }\right)\Big)+&\nonumber\\
  &\dfrac{2}{3}\Big(-2t^*g\left(w_{n},p_{n},q_{n}\right)+t^*g\left(w_{n-1},p_{n-1},q_{n-1}\right)\Big)-E^w_ t,&\\
  &0< \rho<1,\ \ 0< t\leq T,&\nonumber\\
&\dfrac{\partial w_{n+1}}{\partial\rho }\left(0 \right){\rm =0, }~w_{n+1}\left({\rm 1}  \right){\rm =}0,~~~~~0< t\leq T,&\nonumber\\ 
&w_{0}\left(\rho\right){\rm =}0,~~~~~0\leq \rho\leq1,&\nonumber
\end{eqnarray} 
where 
 \begin{equation}\label{ak0llll}
 a'_{n}=0,~~~ a'_k=\dfrac{2( t^*)^{1-\alpha}\Big((k+1)^{1-\alpha}-(k)^{1-\alpha}\Big)}{3\Gamma(2-\alpha)},~~~a'_{k+1}\leq a'_k.
  \end{equation}
   From \eqref{ak0}, one can conclude that there exists positive $C^*_1$ such that
\begin{equation}\label{imp111}
\max\{\|E_ t^c\|_{\infty}, \|E^w_ t\|_{\infty}\}<C_1^*( t^*)^{3-\alpha}.
\end{equation}
By defining  $\xi(\rho, t)$  as follows  
\begin{eqnarray}\label{HM9}
&\frac{\textup{d}}{\textup{d} t }\xi(\rho_0, t)=\nu(\xi(\rho_0, t) , t),~~~0< t\leq T, ~~~\xi \left(\rho_0,0\right)=\rho_0,~~~0\leq\rho\leq1,&
\end{eqnarray}
it is easy to conclude that for each $\rho\in[0,1]$ there exists $\rho_0\in[0,1]$ such that
$\xi(\rho_0, t)=\rho.$  By substituting \eqref{HM9} in \eqref{pqd}--\eqref{YAS2}, we conclude
\begin{equation}\label{YABsser}
\dfrac{\textup{d} p(\xi(\rho_0, t), t)}{\textup{d}  t}=(g_{11}\left(c,w,p,q,d\right)p+g_{12}\left(c,w,p,q,d\right)q +g_{13}\left(c,w,p,q,d\right)d)\Big|_{\rho=\xi(\rho_0, t)},
\end{equation}
\begin{equation}
\dfrac{\textup{d} q(\xi(\rho_0, t), t)}{\textup{d}  t}=(g_{21}\left(c,w,p,q,d\right)p+g_{22}\left(c,w,p,q,d\right)q +g_{23}\left(c,w,p,q,d\right)d)\Big|_{\rho=\xi(\rho_0, t)},
\end{equation}
\begin{equation}\label{HYHMsser}
\dfrac{ \textup{d}d(\xi(\rho_0, t), t)}{\textup{d}  t}=(g_{31}\left(c,w,p,q,d\right)p+g_{32}\left(c,w,p,q,d\right)q +g_{33}\left(c,w,p,q,d\right)d)\Big|_{\rho=\xi(\rho_0, t)},
\end{equation}
\[0\le \rho\le 1,\ \ 0< t\leq T,\]
\begin{equation*}
~~~~~p\left(\rho,0\right)=\ p_0\left(\rho\right),~~~q\left(\rho,0\right)=\ q_0\left(\rho\right),~~~ d\left(\rho,0\right)=d_0\left(\rho\right),~~~~{\rm 0}\le \rho \le {\rm 1.\ }
\end{equation*}
By considering  $\rho=\xi(\rho_0, t_{n+1})$, the problem \eqref{YABsser}--\eqref{HYHMsser}, using the Midpoint rule,  becomes
\[  p(\rho, t_{n+1})=p(\smallunderbrace{\rho- 2t^*\nu(\smallunderbrace{\rho- t^*\nu(\rho, t_n)}_{\textup{Forward-difference}},t_n)}_{\textup{Midpont rule}}, t_{n-1})+\]
\[ 2t^*(g_{11}\left(c_{n},w_{n},p_n,q_n,d_n\right)p_n+g_{12}\left(c_{n},w_{n},p_n,q_n,d_n\right)q_n +\]
\begin{equation}\label{YABss}
g_{13}(c_{n},w_{n},p_n,q_n,d_n)d_n)\Big|_{(\rho- t^*\nu(\rho, t_n), t_n)}+E_{ t}^p,
\end{equation}
\[q(\rho, t_{n+1})=q(\rho- 2t^*\nu(\rho- t^*\nu(\rho, t_n),t_n), t_{n-1})+\] \[ 2t^*(g_{21}\left(c_{n},w_{n},p_n,q_n,d_n\right)p_n+g_{22}\left(c_{n},w_{n},p_n,q_n,d_n\right)q_n+\]
\begin{equation}
g_{23}\left(c_{n},w_{n},p_n,q_n,d_n\right)d_n)\Big|_{(\rho- t^*\nu(\rho, t_n), t_n)}+E_{ t}^q,
\end{equation}
\[ d(\rho, t_{n+1})=d(\rho- 2t^*\nu(\rho- t^*\nu(\rho, t_n),t_n), t_{n-1}), t_{n-1})+\] 
\[ 2t^*(g_{31}\left(c_{n},w_{n},p_n,q_n,d_n\right)p_n+g_{32}\left(c_{n},w_{n},p_n,q_n,d_n\right)q_n+\]
\begin{equation}\label{HYHMss}
g_{33}\left(c_{n},w_{n},p_n,q_n,d_n\right)d_n)\Big|_{(\rho- t^*\nu(\rho, t_n), t_n)}+E_{ t}^d,
\end{equation}
\begin{equation}\label{YABww}
~~~~~p\left(\rho,0\right)=\ p_0\left(\rho\right),~~~q\left(\rho,0\right)=\ q_0\left(\rho\right),~~~ d\left(\rho,0\right)=d_0\left(\rho\right),~~~~{\rm 0}\le \rho \le {\rm 1,\ }
\end{equation}
\begin{equation}
\dfrac{1}{\rho^2}\dfrac{\partial }{\partial \rho}\left(\rho^2\dfrac{{v_n}}{R_n}\right)=h\left(c_n,w_n,p_n,q_n,d_n\right),~~v_n(0)=0,\ \ 0<\rho\leq 1,\ \ t>0, \nonumber
\end{equation}
where  $\rho- 2t^*\nu(\rho- t^*\nu(\rho, t_n),t_n)$ is an approximation of $\xi(\rho_0,t_{n-1})$ using the Midpoint rule and $\rho- t^*\nu(\rho, t_n)$  is an approximation of $\xi(\rho_0,t_{n})$ using the forward-difference formula.
Moreover, from \eqref{HM345} we have
 \begin{eqnarray}\label{TM6}  
&\dfrac{dR(t)}{dt}=R(t)\int_0^1\rho^2h(c,w,p,q,d)d\rho,\ \ t>0, & 
\end{eqnarray}
therefore
 \begin{eqnarray}\label{TM6KK}  
&R_{n+1}=R_{n-1}e^{\int_{t_{n-1}}^{t_{n+1}}\int_0^1\rho^2h(c,w,p,q,d)d\rho dt}=R_{n-1}e^{\overbrace{2t^*\int_0^1\rho^2h(c_n,w_n,p_n,q_n,d_n)d\rho}^{\textup{Midpoint rule}} }+E_t^R. & 
\end{eqnarray}
 \eqref{YABss}--\eqref{HYHMss} is derived from \eqref{YABsser}--\eqref{HYHMsser} employing Midpoint rule. Also    $\rho- 2t^*\nu(\rho- t^*\nu(\rho, t_n),t_n)$ is an approximation of $\xi(\rho_0,t_{n-1})$ using the Midpoint rule and $\rho- t^*\nu(\rho, t_n)$  is an approximation of $\xi(\rho_0,t_{n})$ using the forward-difference formula. Therefore,  we conclude that there exists positive constant $C_2^*$ such that
\begin{equation}\label{imp222}
\max\{\|E^p_ t\|_{\infty},\|E^q_ t\|_{\infty},\|E^d_ t\|_{\infty}\}<C_2^*( t^*)^3.
\end{equation}
From \eqref{imp111}, \eqref{TM6KK} and \eqref{imp222}, we deduce that there exists positive constant $C^*$ such that
\begin{equation}\label{imp444}
\max\{ \|E^p_ t\|_{\infty},~\|E_ t^q\|_{\infty}, ~\|E^d_ t\|_{\infty}\}<C^*( t^*)^3,~ \|E^R_ t\|_{\infty}<C^*( t^*)^3,~\max\{\|E_ t^c\|_{\infty},  \|E^w_ t\|_{\infty}\}<C^*( t^*)^{3-\alpha}.
\end{equation}
Now, we  approximate the solution  of  the problem \eqref{YAS}--\eqref{HM345} by $(c^{ap}_{n+1},w^{ap}_{n+1},p^{ap}_{n+1},q^{ap}_{n+1},d^{ap}_{n+1})$, which is the approximated solution of the following problem
\begin{eqnarray}\label{ap11}
 & \mathcal{C}_{n+1}-\dfrac{a'_0D_{1}}{(R^{ap}_{n+1})^2}\dfrac{1}{{\rho }^{2}}\dfrac{\partial }{\partial\rho }\left({\rho }^{2}\dfrac{\partial \mathcal{C}_{n+1}}{\partial\rho }\right)-t^*\dfrac{2(2{v}^{ap}(1,t_{n})-{v}^{ap}(1,t_{n-1}))\rho}{3R^{ap}_{n+1}}\dfrac{\partial \mathcal{C}_{n+1}}{\partial\rho}=&\nonumber\\
  &c^{ap}_{n}-\dfrac{c^{ap}_{n-1}-c^{ap}_n}{3}-\displaystyle\sum_{k=0}^{n-1} (a'_k-a'_{k+1})\Big(\dfrac{D_{1}}{(R^{ap}_{n-k})^2}\dfrac{1}{{\rho }^{2}}\dfrac{\partial }{\partial\rho }\left({\rho }^{2}\dfrac{\partial c^{ap}_{n-k}}{\partial\rho }\right)\Big)+&\nonumber\\
  &\dfrac{2}{3}\Big(-2t^*f\left(c^{ap}_{n},p^{ap}_{n},q^{ap}_{n}\right)+t^*f\left(c^{ap}_{n-1},p^{ap}_{n-1},q^{ap}_{n-1}\right)\Big),&\\
  &0< \rho<1,\ \ 0< t\leq T,&\nonumber\\
&\dfrac{\partial \mathcal{C}_{n+1}}{\partial\rho }\left(0 \right){\rm =0, }~\mathcal{C}_{n+1}\left({\rm 1}  \right){\rm =}0,~~~~~~\ 0< t\leq T,&\nonumber\\
&\mathcal{C}_{0}\left(\rho \right){\rm =}0,~~~~~0\leq \rho\leq1,&\nonumber
\end{eqnarray}
\begin{eqnarray}\label{ap22}
 & \mathcal{W}_{n+1}-\dfrac{a'_0D_{2}}{(R^{ap}_{n+1})^2}\dfrac{1}{{\rho }^{2}}\dfrac{\partial }{\partial\rho }\left({\rho }^{2}\dfrac{\partial \mathcal{W}_{n+1}}{\partial\rho }\right)-t^*\dfrac{2(2{v}^{ap}(1,t_{n})-{v}^{ap}(1,t_{n-1}))\rho}{3R^{ap}_{n+1}}\dfrac{\partial \mathcal{W}_{n+1}}{\partial\rho}=&\nonumber\\
  &w^{ap}_{n}-\dfrac{w^{ap}_{n-1}-w^{ap}_n}{3}-\displaystyle\sum_{k=0}^{n-1} (a'_k-a'_{k+1})\Big(\dfrac{D_{2}}{(R^{ap}_{n-k})^2}\dfrac{1}{{\rho }^{2}}\dfrac{\partial }{\partial\rho }\left({\rho }^{2}\dfrac{\partial w^{ap}_{n-k}}{\partial\rho }\right)\Big)+&\nonumber\\
 &\dfrac{2}{3}\Big(-2t^*g\left(w^{ap}_{n},p^{ap}_{n},q^{ap}_{n}\right)+t^*g\left(w^{ap}_{n-1},p^{ap}_{n-1},q^{ap}_{n-1}\right)\Big),&\\
  &0< \rho<1,\ \ 0< t\leq T,&\nonumber\\
&\dfrac{\partial \mathcal{W}_{n+1}}{\partial\rho }\left(0 \right){\rm =0, }~\mathcal{W}_{n+1}\left({\rm 1}  \right){\rm =}0,~~~~~0< t\leq T,&\nonumber\\ 
&\mathcal{W}_{0}\left(\rho \right){\rm =}0,~~~~~0\leq \rho\leq1,&\nonumber
\end{eqnarray}
\begin{equation}\label{ap33}
p^{ap}(\rho, t_{n+1})=p^{ap}(\rho- 2t^*\nu_n^{ap}(\rho- t^*\nu^{ap}_n(\rho)), t_{n-1})+
\end{equation}
\[2t^*\Big(g_{11}\left(c^{ap}_{n},w^{ap}_{n},p^{ap}_{n},q^{ap}_n,d^{ap}_{n}\right)p^{ap}_n+g_{12}\left(c^{ap}_{n},w^{ap}_{n},p^{ap}_{n},q^{ap}_n,d^{ap}_{n}\right)q^{ap}_n+\]
\[ g_{13}\left(c^{ap}_{n},w^{ap}_{n},p^{ap}_{n},q^{ap}_n,d^{ap}_{n}\right)d^{ap}_n\Big)\Big|_{(\rho-t^*\nu^{ap}_n(\rho), t_n)},\]  
\begin{equation}\label{ap44}
q^{ap}(\rho, t_{n+1})=q^{ap}(\rho- 2t^*\nu_n^{ap}(\rho- t^*\nu^{ap}_n(\rho)), t_{n-1})+
\end{equation}
\[2t^*\Big(g_{21}\left(c^{ap}_{n},w^{ap}_{n},p^{ap}_{n},q^{ap}_n,d^{ap}_{n}\right)p^{ap}_n+g_{22}\left(c^{ap}_{n},w^{ap}_{n},p^{ap}_{n},q^{ap}_n,d^{ap}_{n}\right)q^{ap}_n+\]
\[ g_{23}\left(c^{ap}_{n},w^{ap}_{n},p^{ap}_{n},q^{ap}_n,d^{ap}_{n}\right)d^{ap}_n\Big)\Big|_{(\rho-t^*\nu^{ap}_n(\rho), t_n)},\] 
\begin{equation}\label{ap55} 
d^{ap}(\rho, t_{n+1})=d^{ap}(\rho- 2t^*\nu_n^{ap}(\rho- t^*\nu^{ap}_n(\rho)), t_{n-1})+
\end{equation}
\[2t^*\Big(g_{31}\left(c^{ap}_{n},w^{ap}_{n},p^{ap}_{n},q^{ap}_n,d^{ap}_{n}\right)p^{ap}_n+g_{32}\left(c^{ap}_{n},w^{ap}_{n},p^{ap}_{n},q^{ap}_n,d^{ap}_{n}\right)q^{ap}_n+\]
\[ g_{33}\left(c^{ap}_{n},w^{ap}_{n},p^{ap}_{n},q^{ap}_n,d^{ap}_{n}\right)d^{ap}_n\Big)\Big|_{(\rho-t^*\nu^{ap}_n(\rho), t_n)},\] 
\begin{equation}\label{ini}
p^{ap}_0=p_0,~~~q^{ap}_0=q_0,~~~d^{ap}_0=d_0,
\end{equation}
\begin{equation}
\dfrac{1}{\rho^2}\dfrac{\partial }{\partial \rho}\left(\rho^2\dfrac{{v}_{n}^{ap}}{R_{n}^{ap}}\right)=h\left(c^{ap}_n,w^{ap}_n,p^{ap}_n,q^{ap}_n,d^{ap}_n\right),~~v^{ap}_n(0)=0,\ \ 0<\rho\leq 1,\ \ t>0, \nonumber
\end{equation}
 \begin{eqnarray}\label{TM6KK**}  
&R^{ap}_{n+1}=R^{ap}_{n-1}e^{\overbrace{2t^*\int_0^1\rho^2h(c^{ap}_n,w^{ap}_n,p^{ap}_n,q^{ap}_n,d^{ap}_n)d\rho}^{\textup{Midpoint rule}} },\ \ t>0, & 
\end{eqnarray}
where $(c_{n+1}^{ap},w_{n+1}^{ap})$ is obtained as an approximation of $(\mathcal{C}_{n+1},\mathcal{W}_{n+1})$ by solving \eqref{ap11}--\eqref{ap22} employing the spectral method.
\noindent  In the following, it is  assumed that
\[(f_1,f_2)_{\omega^{\alpha,\beta}}=\int_{0}^1 \rho^\alpha(1-\rho)^\beta  f_1(\rho, t)f_2(\rho, t)d\rho,~~\|f\|^2_{\omega^{\alpha,\beta}}=(f,f)_{\omega^{\alpha,\beta}.}\]
 We approximate $\mathcal{C}_{n+1}$ employing $\{P_{j}(\rho)\}_{j=0}^{\infty}$, which are chosen such that for each $k\in\mathbb{N}_0$,
 \begin{equation}\label{fb1}
{span}\{P_{0}(\rho),P_{1}(\rho),\cdots,P_{k}(\rho)\}=\Big\{u\in span\{1,\rho,\cdots,\rho^{k+2}\}\mid \dfrac{\partial u(\rho)}{\partial\rho}|_{\rho=0}=0,~u(1)=0\Big\},
\end{equation}
then, we consider $c_{n+1}^{ap}=c^{N^1}_{n+1}$ as follows
\begin{equation}\label{HM6}
c^{N^1}_{n+1}(\rho)=\sum_{j=0}^{N^1}a^{n+1}_jP_j(\rho).
\end{equation}
We  calculate $\{a^{n+1}_j\}_{j=0}^{N^1}$  using the spectral method from 
\begin{equation}\label{collo2}
\Pi_{N^1}^{0,0}c^{N^1}_{n+1}-\Pi_{N^1}^{0,0}\Big(\dfrac{a'_0D_{1}}{(R^{ap}_{n+1})^2}\dfrac{1}{{\rho }^{2}}\dfrac{\partial }{\partial\rho }\Big({\rho }^{2}\dfrac{\partial c^{N^1}_{n+1}}{\partial\rho }\Big)\Big)-t^*\Pi_{N^1}^{0,0}\dfrac{2(2{v}^{ap}(1,t_{n})-{v}^{ap}(1,t_{n-1}))\rho}{3R^{ap}_{n+1}}\dfrac{\partial c^{N^1}_{n+1}}{\partial\rho} =I_{N^1}^{0,0}g^*_n, 
\end{equation}
where $\Pi_{N^1}^{0,0}$ and $I_{N^1}^{0,0}$ are the orthogonal projection and  Jacobi-Gauss-Lobatto interpolation operator with respect to $\rho$ on $[0,1]$, respectively, and  
\[{g^*_n}=\Pi_{N^1}^{0,0}(c^{ap}_{n}-\dfrac{c^{ap}_{n-1}-c^{ap}_n}{3})-\displaystyle\sum_{k=0}^{n-1} (a'_k-a'_{k+1})\Big(\dfrac{D_{1}}{(R^{ap}_{n-k})^2}\dfrac{1}{{\rho }^{2}}\dfrac{\partial }{\partial\rho }\left({\rho }^{2}\dfrac{\partial c^{ap}_{n-k}}{\partial\rho }\right)\Big)+\]
\begin{eqnarray}\label{ap00}
\smallunderbrace{\dfrac{2}{3}\Big(-2t^*f\left(c^{ap}_{n},p^{ap}_{n},q^{ap}_{n}\right)+t^*f\left(c^{ap}_{n-1},p^{ap}_{n-1},q^{ap}_{n-1}\right)\Big)}_{f_n^*}.&
\end{eqnarray}
 Also, we calculate $w^{N^1}_{n+1}$  similar to $c^{N^1}_{n+1}$ to approximate $\mathcal{W}_{n+1}$, the solution of \eqref{ap22}.
Now, using  the principle of mathematical induction, we want to show that for  $k=0,1,\cdots,M,$ there exists positive constants  $c^*$, $w^*$, $p^*$, $q^*$, $d^*$ and $R^*$ such that 
 \begin{equation*}
| p^{ap}_k-p_k|<p^*,~|q^{ap}_k-q_k|< q^*,~|d^{ap}_k-d_k|<d^*,~~p^{ap}_k+q^{ap}_k+d^{ap}_k=N,
\end{equation*}
and
 \begin{equation*}
|c^{ap}_k-c_k|\leq {c^*},~~~|w^{ap}_k-w_k|\leq {w^*},~|R^{ap}_k-R_k|\leq R^*,
\end{equation*}
 where $c$, $w$, $p$, $q$, $d$ and $R$ are the exact solutions of \eqref{YAS}--\eqref{HM345}, respectively. First, we assume that 
 \begin{equation}\label{A1}
| p^{ap}_k-p_k|<p^*,~|q^{ap}_k-q_k|< q^*,~|d^{ap}_k-d_k|<d^*,~~p^{ap}_k+q^{ap}_k+d^{ap}_k=N, ~~~k\leq n<M,
\end{equation}
and
 \begin{equation}\label{A1**}
|c^{ap}_k-c_k|\leq {c^*},~~~|w^{ap}_k-w_k|\leq {w^*},~|R^{ap}_k-R_k|\leq R^*, ~~~k\leq n<M.
\end{equation}
\begin{lemma}\label{convpartial}
Let $c$ be the exact solution of \eqref{YAS} on $[0,1]\times[0,T]$, $c_{n+1}^{ap}=c_{n+1}^{N^1}$ and $\dfrac{\partial^2 c}{\partial\rho^2}$  be a $C^1$-smooth function, also for each $1\leq k\leq n$ the conditions presented in \eqref{A1}--\eqref{A1**} be satisfied. Then, there exists a  positive constant  $K_4$ such that 
\[\|\dfrac{\partial(c^{N^1}_{n+1}- c_{n+1,1}^{N^1})}{\partial\rho}\|^2_{\omega ^{0,0}}\leq \|\dfrac{\partial(c^{ap}_{n}- c_{n,1}^{N^1})}{\partial\rho}\|^2_{\omega ^{0,0}}+\]
\[ {K_4}\Big(({t^*})^{1+\alpha}\sum_{k=0}^n(\dfrac{1}{3})^{\frac{k}{2}}\Big(\|{c_{k,1}^{N^1}}-{c^{ap}_k }\|^2_{\omega^{0,0}}+\|{{p_k }}-{{p^{ap}_k }}\|^2_{\omega^{0,0}}+\|q_k -q^{ap}_k\|^2_{\omega^{0,0}}+\|d_k -d^{ap}_k\|^2_{\omega^{0,0}}+\|R_k -R^{ap}_k\|^2_{\omega^{0,0}}\Big)+\]
\begin{equation}\label{HM1}
\|e_t^c\|^2_{\omega^{0,0}}+t^*K_1^*(N^1)\Big),
\end{equation}
where 
\[\lim_{N^1\rightarrow\infty}K_1^*(N^1)=0,~~\|e_t^c\|_{\infty}\leq (t^*)^{2+\frac{1-\alpha}{2}},\]
and $c_1^{N^1}$  is a polynomial such that
\[c^{N^1}_{n,1}=c_1^{N^1}(\rho,t_n),~~I_{N^1}^{0,0}c^{N^1}_1=c_1^{N^1},~\dfrac{\partial c_1^{N^1}}{\partial\rho }\left(0, t \right){\rm =0, }~c_1^{N^1}\left({\rm 1,} t \right){\rm =}0,~c_1^{N^1}\left(\rho ,0\right){\rm =}0, ~~0\leq t\leq T,\]
and
\begin{equation*}
\lim_{N^1\rightarrow\infty} \Big(\|\dfrac{\partial c_1^{N^1}}{\partial\rho}-\dfrac{\partial c}{\partial\rho}\|^2_{\omega ^{0,0}}+\|I_{N^1}^{0,0}c-c_1^{N^1}\|^2_{\omega^{0,0}}+\|I_{N^1}^{0,0}\dfrac{\partial(c-c_1^{N^1})}{\rho\partial\rho}\|^2_{\omega^{0,0}}+\|I_{N^1}^{0,0}\dfrac{\partial^2(c-c_1^{N^1})}{\partial\rho^2}\|^2_{\omega^{0,0}}\Big)=0.
\end{equation*}
Moreover, if $\dfrac{\partial^2 c}{\partial\rho^2}$ is a $C^m$-smooth function  then
 \[\|\dfrac{\partial(c^{ap}_{n+1}- c_{n+1,1}^{N^1})}{\partial\rho}\|^2_{\omega ^{0,0}}\leq \|\dfrac{\partial(c^{ap}_{n}- c_{n,1}^{N^1})}{\partial\rho}\|^2_{\omega ^{0,0}}+ \]
\[K_4\Big(({t^*})^{1+\alpha}\sum_{k=0}^n(\dfrac{1}{3})^{\frac{k}{2}}\Big(\|{c_{k,1}^{N^1}}-{c^{ap}_k }\|^2_{\omega^{0,0}}+\|{{p_k }}-{{p^{ap}_k }}\|^2_{\omega^{0,0}}+\|q_k -q^{ap}_k\|^2_{\omega^{0,0}}+\|d_k -d^{ap}_k\|^2_{\omega^{0,0}}+\|R_k -R^{ap}_k\|^2_{\omega^{0,0}}\Big)+\]
\begin{equation}\label{HM1**KF}
\|e_t^c\|^2_{\omega^{0,0}}+t^*\smallunderbrace{\dfrac{1}{({N^1})^{2m}}}_{K^*_1(N^1)}\Big).
\end{equation}
\end{lemma}
{\it{Proof}} See Appendix.
\qed

\noindent Similar to the proof of Lemma \ref{convpartial}, we can show that there exists positive constant $K_5$ such that 
\[\|\dfrac{\partial(w^{ap}_{n+1}- w_{n+1,1}^{N^1})}{\partial\rho}\|^2_{\omega ^{0,0}}\leq \|\dfrac{\partial(w^{ap}_{n}- w_{n,1}^{N^1})}{\partial\rho}\|^2_{\omega ^{0,0}}+\]
\[ {K_5}\Big(({t^*})^{1+\alpha}\sum_{k=0}^n(\dfrac{1}{3})^{\frac{k}{2}}\Big(\|{w_{k,1}^{N^1}}-{w^{ap}_k }\|^2_{\omega^{0,0}}+\|{{p_k }}-{{p^{ap}_k }}\|^2_{\omega^{0,0}}+\|q_k -q^{ap}_k\|^2_{\omega^{0,0}}+\|d_k -d^{ap}_k\|^2_{\omega^{0,0}}+\|R_k -R^{ap}_k\|^2_{\omega^{0,0}}\Big)+\]
\begin{equation}\label{HM2}
\|e_t^w\|^2_{\omega^{0,0}}+t^*K_2^*(N^1)\Big),
\end{equation}
where  
\[\lim_{N^1\rightarrow\infty}K_2^*(N^1)=0.\]

\begin{theorem}\label{mainthm}
Let  $c^{ap}_{n+1}=c^{N^1}_{n+1}$ and $w^{ap}_{n+1}=w^{N^1}_{n+1}$.  Then, under assumptions of Lemma \ref{convpartial}, there exist positive constants $M$, $C^*_4$ and $C_5^*$ such that
\begin{equation}\label{imp10101}
\max_{k=0,1,\cdots,n+1}\{\mathfrak{E}_{k}\}\leq C_4^* e^{MT}\Big(( t^*)^{2-\frac{\alpha}{2}}+ K^*_4(N^1)\Big),
\end{equation}
and 
\begin{equation}\label{imp1111}
\max_{k=0,1,\cdots,n+1}\{\mathfrak{e}_{k}\}\leq C_5^* e^{MT}\Big(( t^*)^{2-\frac{\alpha}{2}}+ K^*_4(N^1)\Big),
\end{equation}
where
\[\lim_{N^1\rightarrow\infty}K^*_4(N^1)=0,\]
\begin{equation}\label{defpqd1}
\mathfrak{E}_k=\| p^{ap}_k-p_k\|_{\infty}+\|q^{ap}_k-q_k\|_{\infty}+\|d^{ap}_k-d_k\|_{\infty}+|R^{ap}_{k}-R_{k}|,
\end{equation}
and 
\begin{equation}\label{defcw2}
\mathfrak{e}_k=\| \dfrac{\partial(c^{ap}_k-c_k)}{\partial\rho}\|_{\omega^{0,0}}+\|\dfrac{\partial(w^{ap}_k-w_k)}{\partial\rho}\|_{\omega^{0,0}},
\end{equation}
and if  $\dfrac{\partial^2 c}{\partial\rho^2}$ and  $\dfrac{\partial^2 w}{\partial\rho^2}$ are $C^m$-smooth functions with respect to $\rho$, then
\[K^*_4(N^1)=\dfrac{1}{(N^1)^{m}}.\]
\end{theorem}
{\it{Proof}} See Appendix.
\qed
\noindent Employing the General Sobolev inequalities, there exists a positive constant $C_2$ such that  for $0\leq t\leq T$, we have
\begin{equation}\label{Q}
|c^{ap}_{n+1}-c_{n+1}|\leq C_2(\|c^{ap}_{n+1}-c_{n+1}\|_{\omega^{0,0}}+\|\dfrac{\partial(c^{ap}_{n+1}-c_{n+1})}{\partial\rho}\|_{\omega^{0,0}}).
\end{equation}
Now, using the principle of mathematical induction, Theorem \ref{mainthm}, \eqref{A1} and \eqref{A1**}, we conclude that
 \begin{equation}\label{omdG*}
| p^{ap}_k-p_k|<p^*,~|q^{ap}_k-q_k|< q^*,~|d^{ap}_k-d_k|<d^*,~~p^{ap}_k+q^{ap}_k+d^{ap}_k=N, ~~~0\leq k\leq M,
\end{equation}
and
 \begin{equation}\label{omdG}
|c^{ap}_k-c_k|\leq {c^*},~~~|w^{ap}_k-w_k|\leq {w^*},~|R^{ap}_k-R_k|\leq R^*, ~~~0\leq k\leq M.
\end{equation} 
\noindent Thus from Theorem \ref{mainthm}, we can conclude that the sequence $\{(c^{ap}_n,w^{ap}_n,p^{ap}_n,q^{ap}_n,d^{ap}_n,R^{ap}_n)\}_{n=0}^{\infty}$ converges to the exact solution of problem (\ref{YAS})--(\ref{HM345}) on $[0,1]\times[0,T]$. 

\section{Stability} 
In this section, we want to prove the stability of the presented method. For this aim, first we consider the following problem
\begin{equation}\label{YAS88900}
\dfrac{\partial c}{\partial t }=D_{1}\dfrac{\partial^\alpha}{\partial t^\alpha}\left(\dfrac{1}{(R(t)\rho)^2}\dfrac{\partial }{\partial\rho }\Big({\rho }^{2}\dfrac{\partial c}{\partial\rho }\Big)\right)+\dfrac{v(1,t)\rho}{R(t)}\dfrac{\partial c}{\partial\rho} -f\left(c,p,q\right)+f_1(\rho,t),~ {\rm \ 0}<\rho<1,\  t {\rm >}0,
\end{equation}
\begin{equation}
\dfrac{\partial c}{\partial\rho }\left(0, t \right){\rm =0, }~c\left({\rm 1,} t \right){\rm =}\overline{c}( t), ~~\  t {\rm >}0,
\end{equation} 
\begin{equation}
c\left(\rho ,0\right){\rm =}c_0\left(\rho \right),~{\rm  \ 0}\le \rho \le {\rm 1,}
\end{equation}
\begin{equation}\label{YAS*88900}
\dfrac{\partial w}{\partial t }=D_2\dfrac{\partial^\alpha}{\partial t^\alpha}\left(\dfrac{1}{{(R(t)\rho) }^{2}}\dfrac{\partial }{\partial\rho }\Big({\rho }^{2}\dfrac{\partial w}{\partial\rho }\Big)\right)+\dfrac{v(1,t)\rho}{R(t)}\dfrac{\partial w}{\partial\rho}-g\left(w,p,q\right)+f_2(\rho,t),~ {\rm \ 0}<\rho <1,\  t {\rm >}0,
\end{equation}
\begin{equation}
\dfrac{\partial w}{\partial\rho }\left(0, t \right){\rm =0, }~w\left({\rm 1,} t \right){\rm =}\overline{w}( t), ~~\  t {\rm >}0,
\end{equation}
\begin{equation}\label{YAS188900}
w\left(\rho ,0\right){\rm =}w_0\left(\rho \right),~{\rm \ \ 0}\le \rho \le {\rm 1,}
\end{equation}
\begin{eqnarray}\label{pqd88900}
&\dfrac{\partial p}{\partial t }{\rm +}\nu\dfrac{\partial p}{\partial\rho }{\rm =}
(g_{{\rm 11}}\left(c,w,p,q,d\right)p{\rm +}g_{{\rm 12}}\left(c,w,p,q,d\right)q{\rm +}g_{{\rm 13}}\left(c,w,p,q,d\right)d)+f_3(\rho,t),&\\  
&\dfrac{\partial q}{\partial t }{\rm +}\nu\dfrac{\partial q}{\partial\rho }{\rm =} 
(g_{{\rm 21}}\left(c,w,p,q,d\right)p{\rm +}g_{{\rm 22}}\left(c,w,p,q,d\right)q{\rm +}g_{{\rm 23}}\left(c,w,p,q,d\right)d)+f_4(\rho,t),&\\  
&\dfrac{\partial d}{\partial t }{\rm +}\nu\dfrac{\partial d}{\partial\rho }{\rm =}
(g_{{\rm 31}}\left(c,w,p,q,d\right)p{\rm +}g_{{\rm 32}}\left(c,w,p,q,d\right)q{\rm +}g_{{\rm 33}}\left(c,w,p,q,d\right)d)+f_5(\rho,t),&\\ 
&{\rm 0}\le \rho \le {\rm 1,\ } t {\rm >}0,& \nonumber
\end{eqnarray}
\begin{equation}\label{YAS288900}
~~~~~~~~~~~~~~~~~~~p\left(\rho,0\right)=\ p_0\left(\rho\right),~~~q\left(\rho,0\right)=\ q_0\left(\rho\right),~~~ d\left(\rho,0\right)=d_0\left(\rho\right),~~~~{\rm 0}\le \rho \le {\rm 1,\ }
\end{equation}
\begin{eqnarray} \label{HM34588900} 
&\dfrac{1}{\rho^2}\dfrac{\partial }{\partial \rho}\left(\rho^2\dfrac{{v}}{R(t)}\right)=h\left(c,w,p,q,d\right)+f_6(\rho,t),\ \ 0<\rho\leq 1,\ \ t>0,& \nonumber\\
&{v}\left(0,t\right)=0,\ \ t>0,& \\
&\dfrac{dR(t)}{dt}={v}\left(1,t\right),\ \ t>0, \ \ R(0)=R_0.& \nonumber
\end{eqnarray}
In the following theorem, the stability of the proposed method is proved.
 \begin{theorem}\label{stability}
Let $\epsilon_1$ be a positive constant and  $|f_{i}|<\epsilon_1$ $(i=1,\cdots,6).$  Then, under assumptions of Lemma \ref{convpartial}, there exist positive constants $M_2$, $C^*_6$ and $C_7^*$ such that
\begin{equation}\label{imp10103}
\max_{k=0,1,\cdots,n+1}\{\mathfrak{E}^p_{k}\}\leq C_6^* e^{M_2T}\Big(( t^*)^{2-\frac{\alpha}{2}}+  K^*_5(N^1)+\epsilon_1\Big),
\end{equation}
and
 \begin{equation}\label{imp11113}
\max_{k=0,1,\cdots,n+1}\{\mathfrak{e}^p_{k}\}\leq C_7^* e^{M_2T}\Big(( t^*)^{2-\frac{\alpha}{2}}+ K^*_5(N^1)+\epsilon_1\Big),
\end{equation}
where
\[\lim_{N^1\rightarrow\infty}K^*_5(N^1)=0,\]
\begin{equation}\label{defpqd1**}
\mathfrak{E}^p_k=\| p^{ap,p}_k-p_k\|_{\infty}+\|q^{ap,p}_k-q_k\|_{\infty}+\|d^{ap,p}_k-d_k\|_{\infty}+|R^{ap,p}_{k}-R_{k}|,
\end{equation}
and 
\begin{equation}\label{defcw2**}
\mathfrak{e}^p_k=\| \dfrac{\partial(c^{ap,p}_k-c_k)}{\partial\rho}\|_{\omega^{0,0}}+\|\dfrac{\partial(w^{ap,p}_k-w_k)}{\partial\rho}\|_{\omega^{0,0}},
\end{equation}
where $(c^{ap,p},w^{ap,p},p^{ap,p},q^{ap,p},d^{ap,p},R^{ap,p})$ is the approximated solution of the perturbed problem \eqref{YAS88900}-\eqref{HM34588900}  using the presented method and  $(c,w,p,q,d,R)$ is the solution of \eqref{YAS}--\eqref{HM345}.
\end{theorem}
{\it{Proof}} See Appendix.
\qed
\section{Numerical experiments}
\begin{figure}
\centering
\includegraphics[scale=.4]{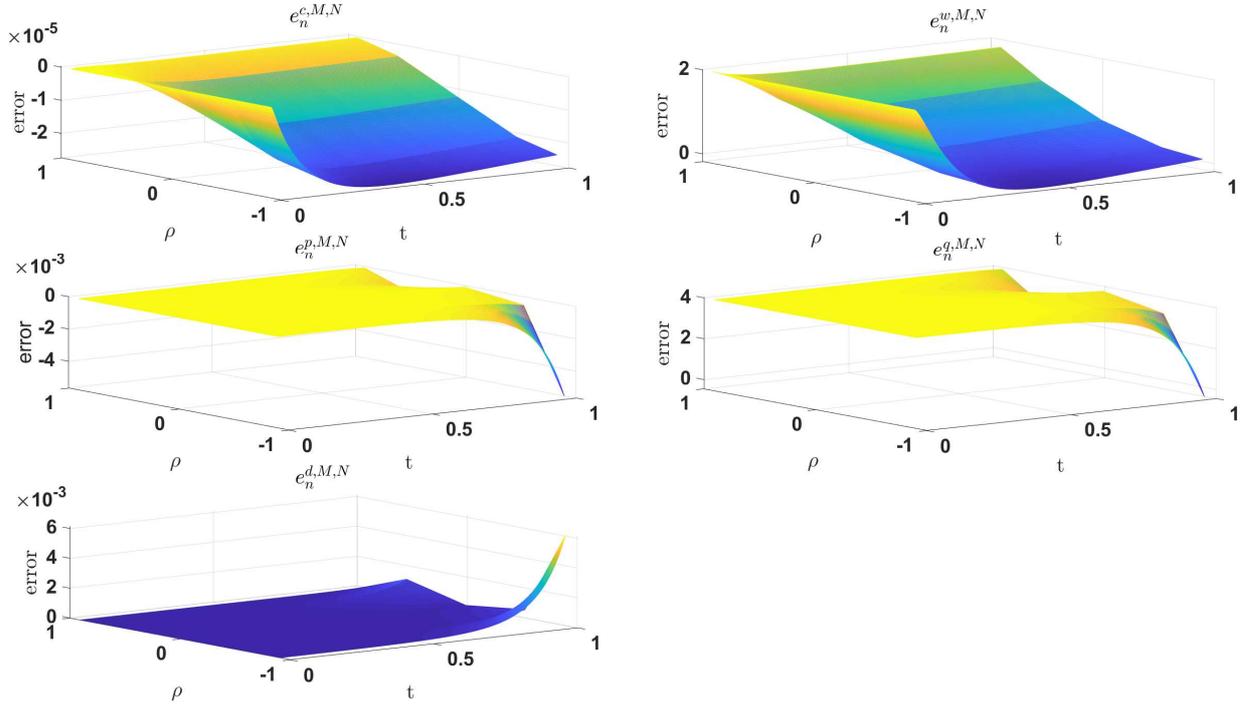}
\caption{Error functions $e_n^{c,M,N}$, $e_n^{w,M,N}$, $e_n^{p,M,N}$, $e_n^{q,M,N}$, $e_n^{d,M,N}$, $e_n^{c,M,N}$ for N=5 and M=200 and $\alpha=0.1$.}\label{exact_err1}
\end{figure}

\begin{figure}
\centering
\includegraphics[scale=.4]{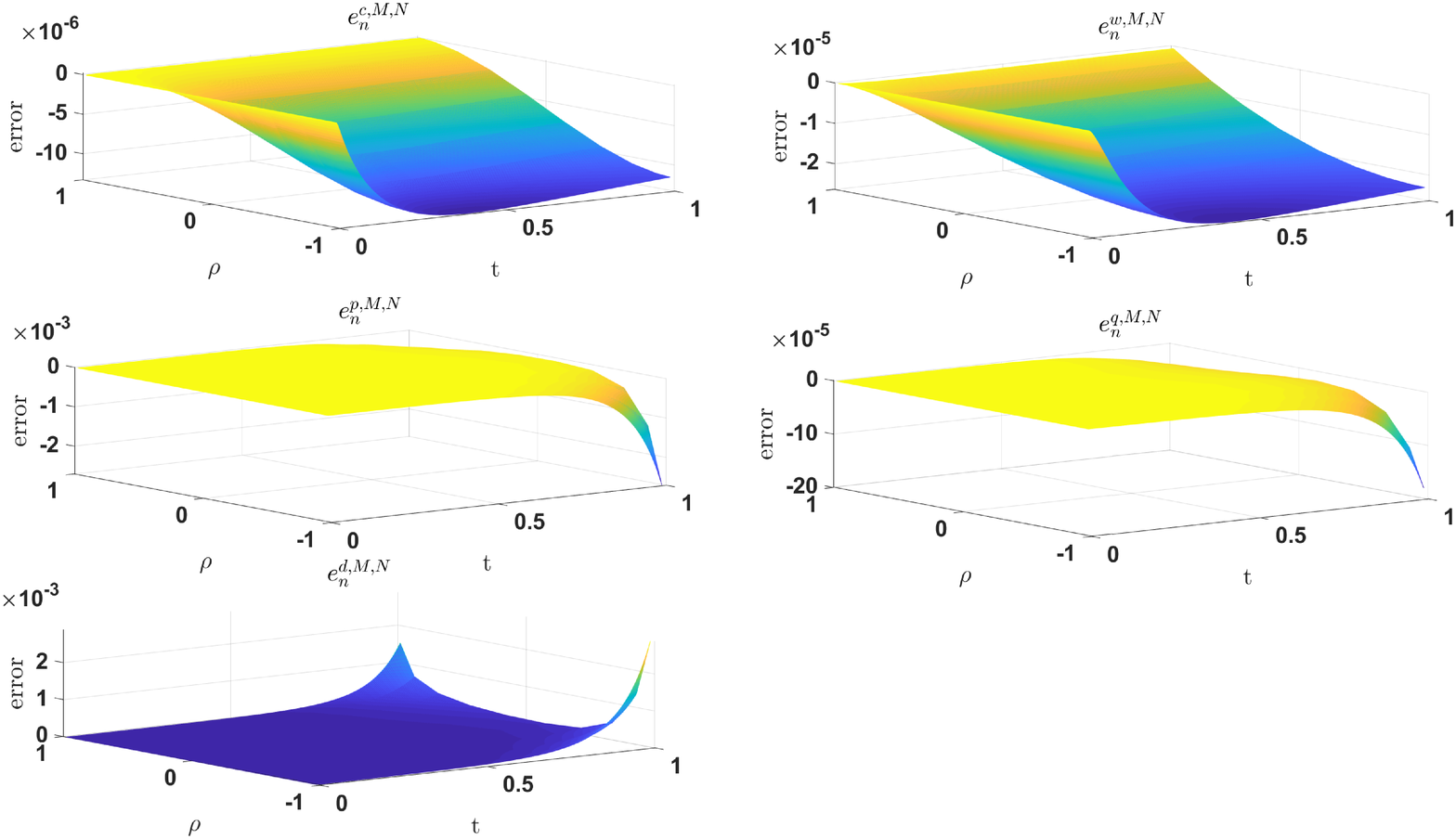}
\caption{Error functions $e_n^{c,M,N}$, $e_n^{w,M,N}$, $e_n^{p,M,N}$, $e_n^{q,M,N}$, $e_n^{d,M,N}$, $e_n^{c,M,N}$ for N=10 and M=300 and $\alpha=0.1$.}\label{exact_err2}
\end{figure}

\begin{figure}
\centering
\includegraphics[scale=.45]{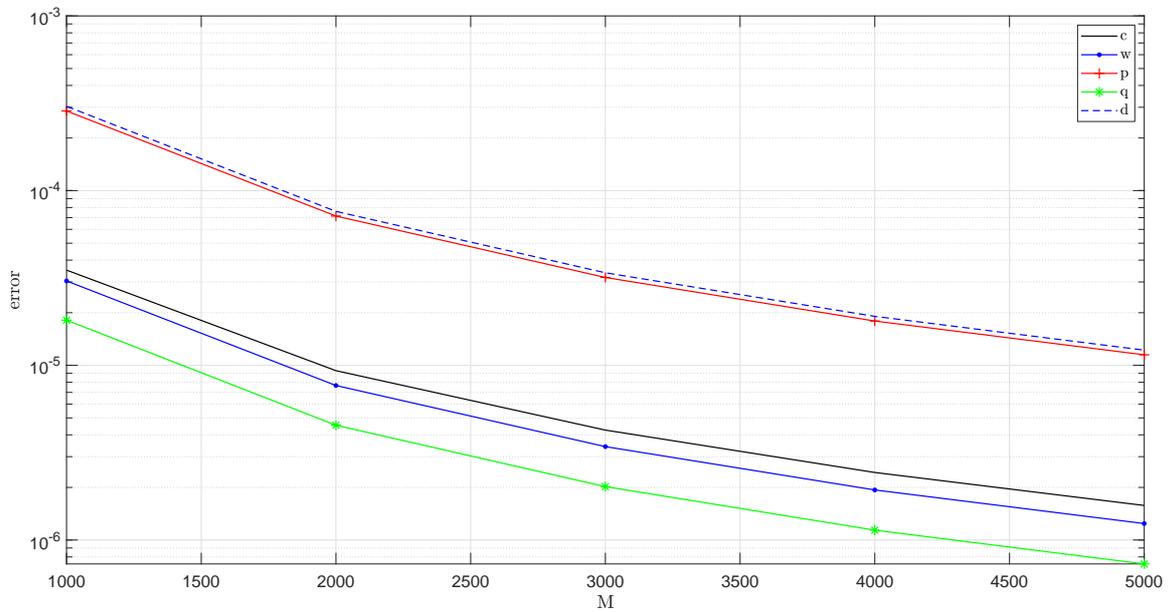}
\caption{Maximum time-errors with N=20 and various M and $\alpha=0.1$.}\label{ertime}
\end{figure}

\begin{figure}
\centering
\includegraphics[scale=.5]{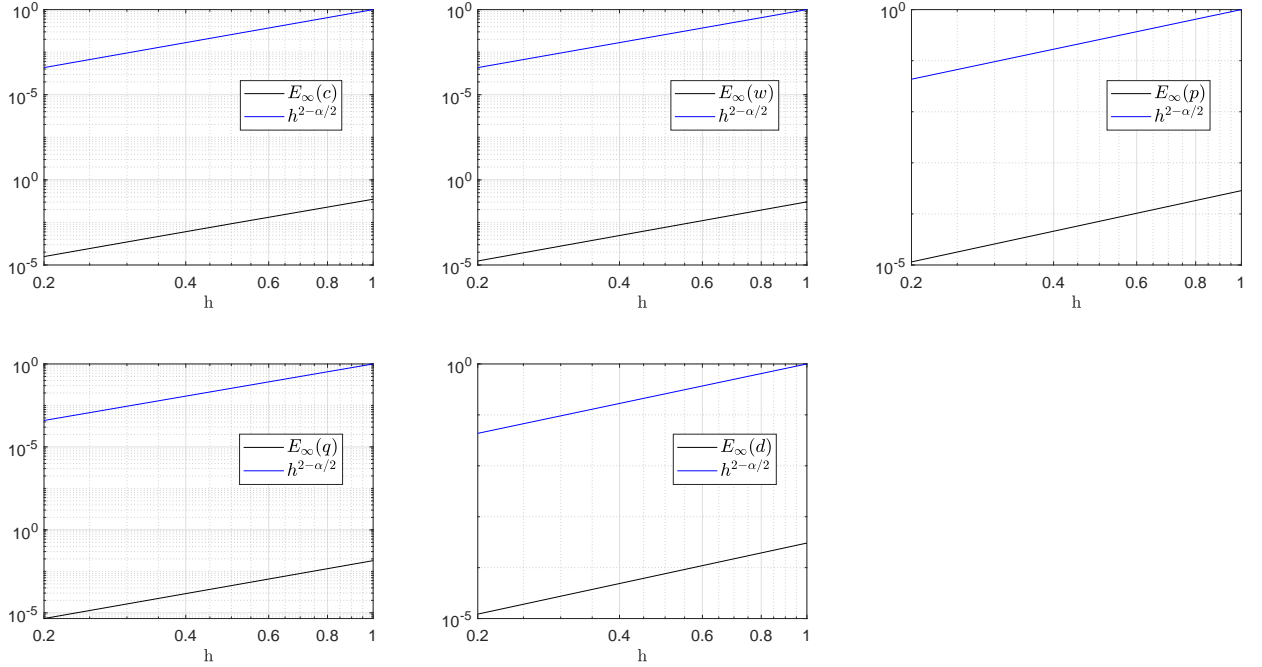}
\caption{The behaviour of time-errors in LogLog scale with $\alpha=0.1$.}
\label{nesbat_orthof}
\end{figure}

\begin{figure}
\centering
\includegraphics[scale=.45]{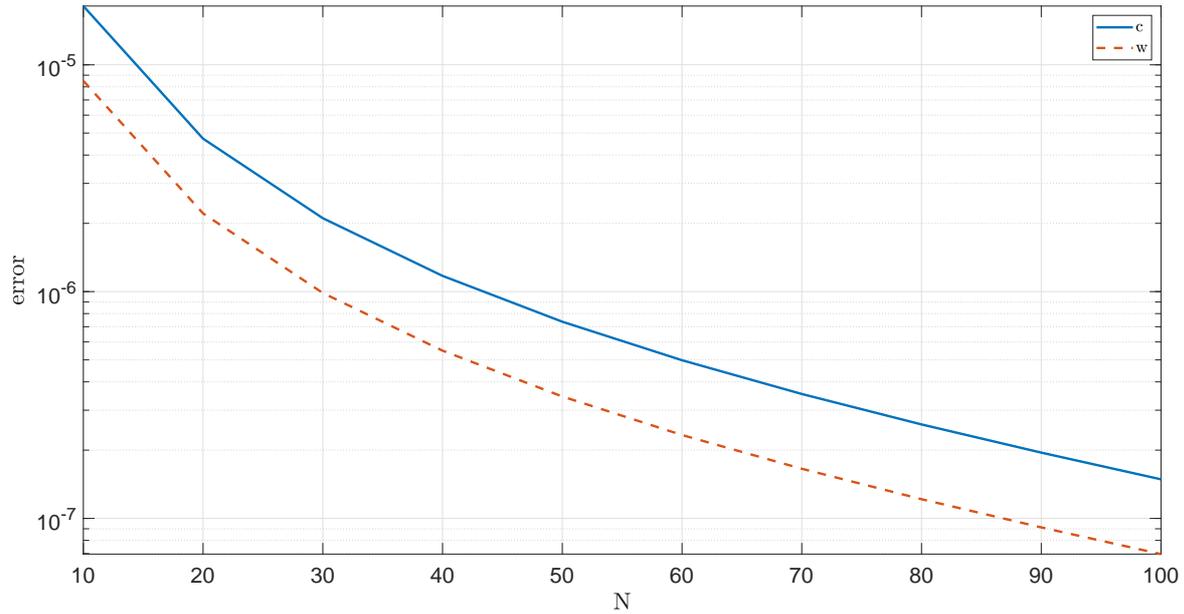}
\caption{Maximum time-error function for M=200 and various N and $\alpha=0.1$.}
\label{erspace}
\end{figure}

In this section, we solve the model of tumor growth by applying finite difference method for approximating the time derivative using mesh points $\{t_i\}_{i=1}^{M}$  where $t_i=i\frac{T}{M}$ and $M$ is a positive integer and spectral method in space.  To construct trial functions for spectral method, which satisfy the boundary conditions, we use orthogonal Legendre polynomials as trial functions in the form of (\ref{fb1}) on (\ref{YAS})-(\ref{YAS1}) as follows
\[c_{n+1}^N(\rho)=\sum_{i=0}^N c_i^{n+1} p_i(\rho) ,~~w_{n+1}^N(\rho)=\sum_{i=0}^N w_i^{n+1} p_i(\rho) ,\]
where
\begin{equation}\label{tf}
p_i(\rho)= L_i (\rho) -\dfrac{2i+3}{(i+2)^2}L_{i+1} (\rho) -\left( \dfrac{i+1}{i+2}\right) ^2L_{i+2} (\rho) , ~~i=0,\ldots, N,
\end{equation}
\textit{Remark: The scaling factors in the trial functions (\ref{tf}) play the role of precondition factor for the collocation matrices \cite{huang2018spectral}}.\\
Also the Gauss quadrature points $\{x_i^{0,0}\}_{i=1}^{N}$ (i.e., the zeros of Legendre polynomial of degree $N+1$) are considered as collocation points.\\
In order to verify our numerical results, we need to present the following definition.
\begin{definition}
A sequence $\{x_n\}_{n=1}^{\infty}$ is said to converge to $x$ with order $p$ if there exists a constant $C$ such that $\vert x_n-x\vert \leq Cn^{-p} ,~~\forall n$.
This can be written as $\vert x_n-x\vert=\mathcal{O}(n^{-p})$.  A practical method to calculate the rate of convergence for a discretization method is to use the following formula
\begin{equation}\label{ratio_remark}
p\approx \dfrac{\log_e(e_{n_2}/e_{n_1})}{\log_e(n_1/n_2)},
\end{equation}
where $e_{n_1}$ and $e_{n_2}$ denote the errors with respect to the step sizes $\dfrac{1}{n_1}$ and $\dfrac{1}{n_2}$, respectively \cite{gautschi1997numerical}.
\end{definition}
Now, using these trial functions, we want to solve the following example.\\
\textbf{Example}1. Consider the following problem:
\begin{equation*}\label{Solve}
\dfrac{\partial c}{\partial t }=\dfrac{\partial^\alpha }{\partial t^\alpha }\left(\dfrac{{1}}{{12}(R(t)\rho)^2}\dfrac{\partial }{\partial\rho }\Big({\rho }^{2}\dfrac{\partial c}{\partial\rho }\Big)\right)+\dfrac{{v}(1,t)\rho}{R(t)}\dfrac{\partial c}{\partial\rho}+\dfrac{c}{16}+\dfrac{12p}{88},~ {\rm \ 0}<\rho<1,\  t {\rm >}0,
\end{equation*}
\begin{equation*}
\dfrac{\partial c}{\partial\rho }\left(0, t \right){\rm =0, }~c\left({\rm 1,} t \right){\rm =}0, ~~\  t {\rm >}0,
\end{equation*} 
\begin{equation}\label{ex1eq1}
c\left(\rho ,0\right){\rm =}0, ~~ 0\leq \rho\leq 1,
\end{equation}
\begin{equation*}
\dfrac{\partial w}{\partial t }=\dfrac{\partial^\alpha }{\partial t^\alpha }\left(\dfrac{1}{12{(R(t)\rho) }^{2}}\dfrac{\partial }{\partial\rho }\Big({\rho }^{2}\dfrac{\partial w}{\partial\rho }\Big)\right)+\dfrac{{v}(1,t)\rho}{R(t)}\dfrac{\partial w}{\partial\rho}+\dfrac{3c}{115}+\dfrac{12p}{188},~ {\rm \ 0}<\rho <1,\  t {\rm >}0,
\end{equation*}
\begin{equation*}
\dfrac{\partial w}{\partial\rho }\left(0, t \right){\rm =0, }~w\left({\rm 1,} t \right){\rm =}0, ~~\  t {\rm >}0,
\end{equation*}
\begin{equation}\label{ex1eq2}
w\left(\rho ,0\right){\rm =}0\left(\rho \right),~{\rm \ \ 0}\le \rho \le {\rm 1,}
\end{equation}
\begin{eqnarray}
&\dfrac{\partial p}{\partial t }{\rm +}{\dfrac{(v-\rho v(1,t))}{R(t)}}\dfrac{\partial p}{\partial\rho }{\rm =}
(q)p+(\dfrac{c}{2})q +(p)d+f_p, &\\ 
&\dfrac{\partial q}{\partial t }{\rm +}\dfrac{{(v-\rho v(1,t))}}{R(t)}\dfrac{\partial q}{\partial\rho }{\rm =} 
p+(2p)q +f_q,
&\\  
&\dfrac{\partial d}{\partial t }{\rm +}\dfrac{{(v-\rho v(1,t))}}{R(t)}\dfrac{\partial d}{\partial\rho }{\rm =}
p+(p)q+f_d,~~{\rm 0}\le \rho \le {\rm 1,\ } t {\rm >}0,& 
\end{eqnarray}
\begin{equation}
~~~~~~~~~~~~~~~~~~~p\left(\rho,0\right)=\ -((2\rho-1)^2+ 1),~~~q\left(\rho,0\right)=0,~~~ d\left(\rho,0\right)=4\rho,~~~~{\rm 0}\le \rho \le {\rm 1,\ }
\end{equation}
\begin{eqnarray}
&\dfrac{1}{\rho^2}\dfrac{\partial }{\partial \rho}\left(\rho^2\dfrac{{v}}{R(t)}\right)=\dfrac{2p-d}{2}+f_v,\ \ 0<\rho\leq 1,\ \ t>0,& \nonumber\\
&{v}\left(0,t\right)=-\dfrac{\epsilon e^{-1}}{(t+1)^2},\ \ t>0,& \\
&\dfrac{dR(t)}{dt}={v}\left(1,t\right),\ \ t>0, \ \ R(0)=0.5,& \label{ex1eqend}
\end{eqnarray}
and the exact solutions are as follows
\begin{eqnarray}
&c(t,\rho)=4t(2\rho + 1)(\rho - 1)^2,~~w(t,\rho)=-8t(\rho^2 - 1),~~p(t,\rho)=-\exp(t)((2\rho - 1)^2 + 1),&\\
&q(t,\rho)=-t((2\rho - 1)^2 - 1),~~d(t,x)=\exp(t)((2\rho - 1)^2 + 1) + t((2\rho - 1)^2 - 1) + 1,&
\end{eqnarray}
It should be noted that in order to use Legendre polynomials, we map the domain of the problem
(\ref{ex1eq1})-(\ref{ex1eqend}) to $[-1, 1]$.
We carried out the numerical computations using the \textbf{MATLAB 2018a} program using a computer with the Intel Core i7 processor (2.90 GHz, 4 physical cores).\\
In Figures \ref{exact_err1}-\ref{exact_err2}, we have plotted the graph of error functions
\begin{equation*}
e_n^{c,N,M}=c_n^{ap}-c,~~e_n^{w,N,M}=w_n^{ap}-w,
\end{equation*}
\begin{equation*}
e_n^{p,N,M}=p_n^{ap}-p,~~e_n^{q,N,M}=q_n^{ap}-q,
\end{equation*}
 We have presented the maximum time-error by considering constant $N=20$ in Table \ref{timeerrorth}. To better see the time-error of numerical results, Figure \ref{ertime} is presented. Also, in Table \ref{nesbat_orthot} and Figure \ref{nesbat_orthof} the computed order of convergence ($p$) using (\ref{ratio_remark}) for numerical finite difference method is shown. It is shown that the finite difference method has almost $\mathcal{O}(h^{2-\alpha/2})$ error, as someone whould expect from convergence Theorem \ref{mainthm}. In Table \ref{spaceerrorth} we have presented the maximum space-errors  by constant $M=200$ and various values of $N$ .To better see the space-error of numerical results, Figure \ref{erspace} is presented.
\begin{table}
\begin{center}
\scriptsize{\begin{tabular}{lllllllllllll}
\hline
{Error of} &  & \multicolumn{1}{c}{M=100} & \multicolumn{1}{c}{} & \multicolumn{1}{c}{M=1000} & \multicolumn{1}{c}{} & \multicolumn{1}{c}{M=2000} & \multicolumn{1}{c}{} & \multicolumn{1}{c}{M=3000} & \multicolumn{1}{c}{} & \multicolumn{1}{c}{M=4000} & \multicolumn{1}{c}{} & \multicolumn{1}{c}{M=5000} \\ 
\cline{1-1}\cline{3-3}\cline{5-5}\cline{7-7}\cline{9-9}\cline{11-11}\cline{13-13}
\multicolumn{1}{c}{c} &  & \multicolumn{1}{c}{2.22198e-03} & \multicolumn{1}{c}{} & \multicolumn{1}{c}{3.50646e-05} & \multicolumn{1}{c}{} & \multicolumn{1}{c}{9.31974e-06} & \multicolumn{1}{c}{} & \multicolumn{1}{c}{4.25763e-06} & \multicolumn{1}{c}{} & \multicolumn{1}{c}{0.24352e-06} & \multicolumn{1}{c}{} & \multicolumn{1}{c}{1.57668e-06} \\ 
\multicolumn{1}{c}{w} &  & \multicolumn{1}{c}{2.85697e-03} & \multicolumn{1}{c}{} & \multicolumn{1}{c}{3.04010e-05} & \multicolumn{1}{c}{} & \multicolumn{1}{c}{7.66240e-06} & \multicolumn{1}{c}{} & \multicolumn{1}{c}{3.42437e-06 } & \multicolumn{1}{c}{} & \multicolumn{1}{c}{1.93429e-06} & \multicolumn{1}{c}{} & \multicolumn{1}{c}{1.24213e-06} \\ 
\multicolumn{1}{c}{p} &  & \multicolumn{1}{c}{2.89789e-02} & \multicolumn{1}{c}{} & \multicolumn{1}{c}{2.85568e-04} & \multicolumn{1}{c}{} & \multicolumn{1}{c}{7.14999e-05} & \multicolumn{1}{c}{} & \multicolumn{1}{c}{3.18273e-05} & \multicolumn{1}{c}{} & \multicolumn{1}{c}{1.79269e-05} & \multicolumn{1}{c}{} & \multicolumn{1}{c}{1.14866e-05} \\ 
\multicolumn{1}{c}{q} &  & \multicolumn{1}{c}{1.84616e-03} & \multicolumn{1}{c}{} & \multicolumn{1}{c}{1.81339e-05} & \multicolumn{1}{c}{} & \multicolumn{1}{c}{4.54002e-06} & \multicolumn{1}{c}{} & \multicolumn{1}{c}{2.02089e-06} & \multicolumn{1}{c}{} & \multicolumn{1}{c}{1.13827e-06} & \multicolumn{1}{c}{} & \multicolumn{1}{c}{7.29339e-07} \\ 
\multicolumn{1}{c}{d} &  & \multicolumn{1}{c}{3.08257e-02} & \multicolumn{1}{c}{} & \multicolumn{1}{c}{3.03702e-04} & \multicolumn{1}{c}{} & \multicolumn{1}{c}{7.60399e-05} & \multicolumn{1}{c}{} & \multicolumn{1}{c}{3.38482e-05} & \multicolumn{1}{c}{} & \multicolumn{1}{c}{1.90652e-05} & \multicolumn{1}{c}{} & \multicolumn{1}{c}{1.22160e-05} \\ 
\hline
\end{tabular}}
\end{center}
\caption{{Maximum time-errors with N=20 and various M and $\alpha=0.1$. }}
\label{timeerrorth}
\end{table}
\begin{table}
\begin{center}
\scriptsize{\begin{tabular}{|lllllllllll|l|}
\hline
Error ratio (p)  of &  & \multicolumn{1}{c}{c} & \multicolumn{1}{c}{} & \multicolumn{1}{c}{w} & \multicolumn{1}{c}{} & \multicolumn{1}{c}{p} & \multicolumn{1}{c}{} & \multicolumn{1}{c}{q} & \multicolumn{1}{c}{} & \multicolumn{1}{c|}{d} & \multicolumn{1}{c|}{$2-\alpha/2$} \\ 
\cline{1-1}\cline{3-3}\cline{5-5}\cline{7-7}\cline{9-9}\cline{11-12}
\multicolumn{1}{|c}{M=1000,2000} &  & \multicolumn{1}{c}{1.9116} & \multicolumn{1}{c}{} & \multicolumn{1}{c}{1.9882} & \multicolumn{1}{c}{} & \multicolumn{1}{c}{1.9979} & \multicolumn{1}{c}{} & \multicolumn{1}{c}{1.9978} & \multicolumn{1}{c}{} & \multicolumn{1}{c|}{1.9978} & \multicolumn{1}{c|}{} \\ 
\multicolumn{1}{|c}{M=2000,3000} &  & \multicolumn{1}{c}{1.9321} & \multicolumn{1}{c}{} & \multicolumn{1}{c}{1.9863} & \multicolumn{1}{c}{} & \multicolumn{1}{c}{1.9962} & \multicolumn{1}{c}{} & \multicolumn{1}{c}{1.9962} & \multicolumn{1}{c}{} & \multicolumn{1}{c|}{1.9961} & \multicolumn{1}{c|}{} \\ 
\multicolumn{1}{|c}{M=3000,4000} &  & \multicolumn{1}{c}{1.9420} & \multicolumn{1}{c}{} & \multicolumn{1}{c}{1.9854} & \multicolumn{1}{c}{} & \multicolumn{1}{c}{1.9953} & \multicolumn{1}{c}{} & \multicolumn{1}{c}{1.9953} & \multicolumn{1}{c}{} & \multicolumn{1}{c|}{1.9953} & \multicolumn{1}{c|}{1.95} \\ 
\multicolumn{1}{|c}{M=4000,5000} &  & \multicolumn{1}{c}{1.9481} & \multicolumn{1}{c}{} & \multicolumn{1}{c}{1.9848} & \multicolumn{1}{c}{} & \multicolumn{1}{c}{1.9947} & \multicolumn{1}{c}{} & \multicolumn{1}{c}{1.9948} & \multicolumn{1}{c}{} & \multicolumn{1}{c|}{1.9947} & \multicolumn{1}{c|}{} \\ 
\multicolumn{1}{|c}{M=5000,6000} &  & \multicolumn{1}{c}{1.9492} & \multicolumn{1}{c}{} & \multicolumn{1}{c}{1.9839} & \multicolumn{1}{c}{} & \multicolumn{1}{c}{1.9944} & \multicolumn{1}{c}{} & \multicolumn{1}{c}{1.9939} & \multicolumn{1}{c}{} & \multicolumn{1}{c|}{1.9943} & \multicolumn{1}{c|}{} \\ 
\hline
\end{tabular}}
\end{center}
\caption{{The rate of convergence with respect to the time variable and $\alpha=0.1$}}
\label{nesbat_orthot}
\end{table}
\begin{table}
\begin{center}
\scriptsize{\begin{tabular}{lllllllllll}
\hline
\multicolumn{1}{c}{Error of} & \multicolumn{1}{c}{} & \multicolumn{1}{c}{N=10} & \multicolumn{1}{c}{} & \multicolumn{1}{c}{N=20} & \multicolumn{1}{c}{} & \multicolumn{1}{c}{N=40} & \multicolumn{1}{c}{} & \multicolumn{1}{c}{N=80} & \multicolumn{1}{c}{} & \multicolumn{1}{c}{N=100} \\ 
\cline{1-1}\cline{3-3}\cline{5-5}\cline{7-7}\cline{9-9}\cline{11-11}
\multicolumn{1}{c}{c} & \multicolumn{1}{c}{} & \multicolumn{1}{c}{1.81574e-05} & \multicolumn{1}{c}{} & \multicolumn{1}{c}{4.72446e-06} & \multicolumn{1}{c}{} & \multicolumn{1}{c}{1.17313e-06} & \multicolumn{1}{c}{} & \multicolumn{1}{c}{2.59655e-07} & \multicolumn{1}{c}{} & \multicolumn{1}{c}{1.48707e-07} \\ 
\multicolumn{1}{c}{w} & \multicolumn{1}{c}{} & \multicolumn{1}{c}{849925e-06} & \multicolumn{1}{c}{} & \multicolumn{1}{c}{2.21145e-06} & \multicolumn{1}{c}{} & \multicolumn{1}{c}{5.49129e-07} & \multicolumn{1}{c}{} & \multicolumn{1}{c}{1.21541e-07} & \multicolumn{1}{c}{} & \multicolumn{1}{c}{6.96080e-08} \\ 
\hline
\end{tabular}}
\end{center}
\caption{{Maximum space-error with M=200 and $\alpha=0.1$. }}
\label{spaceerrorth}
\end{table}
\section{Conclusions}
In this paper, we have considered a free boundary problem modelling the growth of tumor including two reaction-diffusion equations describing the diffusion of nutrient and drug in the tumor and three hyperbolic equations describing the evolution of tumor cells. Since in the real situation  the subdiffusion of nutrient and drug in the tumor can be found, we have changed the reaction-diffusion equations to the fractional  ones to consider this class of anomalous diffusion and deal with a more reliable model of tumor growth.  After that in order to have a clear vision of the dynamic of tumor growth and the effect of nutrient and drug on the tumor growth, we have solved the fractional problem.  Applying a combination of finite difference method and spectral method, the fractional problem is solved. We have also proved the method is unconditionally convergent and stable, which leads us to trust  the obtained solution.   Finally by giving the numerical  results,    the theoretical statements are justified.
\section*{Appendix}
We provide here the proofs of theorems and lemmas together with some essential mathematical concepts, lemmas and theorems,
which have been used in the mathematical analysis throughout the paper. 

\begin{lem}\label{gronwall} (Discrete Gronwall Lemma) \textup{\cite{Qu}} Let $f_0\geq 0$ and  $\{h_n\}$ and $\{g_n\}$ be  non-negative sequences. if the sequence $\{\psi_n\}$ satisfies
\[\psi_0\leq g_0,\] 
\[\psi_n\leq f_0+\sum_{k=0}^{n-1}g_k+\sum_{k=0}^{n-1}h_k\psi_k,~~~n\geq 1,\]
then we have
\[\psi_n\leq (f_0+\sum_{k=0}^{n-1}g_k)e^{\sum_{k=0}^{n-1}h_k},~~~n\geq 1.\]
\end{lem}

\begin{defin}\textup{\cite{CHQZ}}
Let $ \Omega\subset \mathbb{R}^n$ be  open and bounded  and $\omega$ be a positive continuous function on $\Omega$. We can define weighted $L^p$-norms as follows:

\begin{equation}\label{lpnorm}
~{\left\| u \right\|_{p,{\omega }}} := {\Big( {\int\limits_\Omega {{\omega }(X){{\left| {u\left( X \right)} \right|}^p}dX} } \Big)^{\frac{1}{p}}},~~~~1\leq p<+\infty,
\end{equation} 
\[\|u\|_{\infty}:=\textup{ess}\sup_{X\in \Omega}|u(X)|,\]
where $X = \left( {{x_1},{x_2},...,{x_n}} \right) \in {\mathbb{R}^n}$ and the space of measurable  functions on $\Omega$ for which this norm is finite forms a Banach
space, indicated by $L^p_{w}(\Omega)$.
Moreover, for $L^2_{w}(\Omega)$  
the following inner product is defined 
\begin{equation*}
{\left( {u,v} \right)_{{\omega }}} := { {\int\limits_\Omega {{\omega }(X)u\left( X \right)v\left( X \right)dX} } },\;\; \forall u , v \in {L_{{\omega }}^2(\Omega)},
\end{equation*}
and for simplicity we show the norm  with ${\left\| u \right\|_{{\omega }}}$.

\end{defin}

Now for simplicity, for  $\omega(x)\equiv 1$, we use $L^p(\Omega)$ and $\|u\|_{p}$ instead of $L^p_\omega(\Omega)$ and $\|u\|_{p,\omega}$, respectively.
\begin{defin} \textup{\cite{Z. Wu}} 
Let $ \Omega\subset \mathbb{R}^n$ be open and bounded, $1\leq p<+\infty$ and $Q_T=\Omega\times]0,T[$. We can define the space
$$ W_{p}^{m,k}\left( Q_T \right)=\left\{ { u :\partial _X^lu, \partial^r_tu \in L^p\left( Q_T \right),~~0\leq \left| l \right|\leq m, 0\leq r \le k} \right\},~~k\geq1,
$$
endowed with the norm
\[ \left\| u \right\|_{W_{p}^{m,k}}:={\Big( {\sum_{1\leq r\leq k}\|\partial_t^r u\|^p_{p,\omega}+\sum_{0\leq|l|\leq m} {\| {\partial _X^lu} \|_{p,\omega }^p}} } \Big)^{\frac{1}{p}},\]
\[\omega(X,t)=\textup{sgn}(|\max\{\partial_t^1 u : t\in [0,T]\}|)(1-\dfrac{1}{T})+\dfrac{1}{T},\]
 where $\textup{sgn}(x)$ is the sign  function, $l = ( {{l_1},{l_2},...,{l_n}} )$, $\left| l \right| = \sum\limits_{j = 1}^n {{l_j}} $,
and
$\partial _X^lu = \frac{{{\partial ^{{l_1} + {l_2} + ... + {l_n}}}u}}{{\partial x_1^{{l_1}}\partial x_2^{{l_2}} \cdots \partial x_n^{{l_n}}}}.$ Moreover, for $W_{2}^{m,k}(Q_T)$, the following inner product can be defined
\[ ( u ,v)_{W_{2}^{m,k}}:=\iint_{Q_T}{\Big( {\sum_{1\leq r\leq k}\partial_t^r u\partial_t^r v\omega+\sum_{0\leq|l|\leq m} {\partial _X^lu} {\partial _X^lv}}\omega } \Big)dXdt,~~~\forall u,v\in W_{2}^{m,k}\left( Q_T \right).\]

\end{defin}

\begin{defin} \textup{\cite{J.H.}} Let $\Omega\subseteq \mathbb{R}^3$ be an open set, $p>\dfrac{5}{2}$ and $Q_T=\Omega\times]0,T[$, then we define $D_{p}(\Omega)$ the trace space of $W_{p}^{2,1}(Q_T)$ at $t=0$ as follows
\begin{equation*}
 D_{p}(\Omega)=\{\varphi :~\exists u\in W_{p}^{2,1}(Q_T)~s.t.~ u(., 0) = \varphi \}.
\end{equation*}The norm defined on $D_{p}(\Omega)$ is
\begin{equation*}
||\varphi ||_{D_p }  = \inf\{T^{-\frac{1}{p}} ||u||_{W_{p}^{2,1}} :~ u\in W_{p}^{2,1} (Q_T),~  u(.,0)=\varphi\}.
\end{equation*}
\end{defin}

\begin{lem}\label{lem2*}
For any $u \in \mathbb{P}_N$ on $]a,b[$, we have
\begin{equation*}
{\left\| {\partial _t^ku} \right\|_{{\omega ^{\alpha + k,\beta + k}}}} \lesssim{N^k}{\left\| u \right\|_{{\omega ^{\alpha ,\beta }}}}, \ k\geq 1,~~{\left\| {\partial _tu} \right\|_{{\omega ^{\alpha ,\beta }}}} \lesssim N^2{\left\| u \right\|_{{\omega ^{\alpha ,\beta }}}},
\end{equation*}
where $\omega^{\alpha,\beta}(t)=(b-t)^{\alpha}(t-a)^{\beta}$.
\end{lem}
\begin{proof} See \cite{STL}.
\end{proof}

\begin{lem}\label{lem3*}
Suppose $I=]a,b[$, $ \max\{\alpha+1,\beta+1\}<m\leq N+1$ and $x_j ^{\alpha ,\beta }$, $(0 \leq j \leq N )$, are the Jacobi-Gauss-Lobatto quadrature nodes and $w_j ^{\alpha ,\beta }$, $(0 \leq j \leq N )$, are the Jacobi-Gauss-Lobatto weights. The Jacobi-Gauss-Lobatto interpolation operator is denoted by ${I_N^{\alpha ,\beta }}u$. For any measurable function  $u$ such that $\partial_t^i u\in L^2_{\omega^{\alpha+i,\beta+i}}(I)$,  $(i=0,\cdots,m)$, we have for $0\leq l\leq m$, there exists constant $C$ independent of $m$, $N$, and $u$, such that
\begin{equation}\label{57}
{\left\| {\partial _t^l\left( {u - I_N^{\alpha ,\beta }u} \right)} \right\|_{{\omega ^{\alpha + l,\beta + l}}}} \le C(\dfrac{(N-m+1)!}{N!})^{\frac{1}{2}}(N)^{l-\frac{m+1}{2}}{\left\| {\partial _t^mu} \right\|_{{\omega ^{\alpha + m,\beta + m}}}}.
\end{equation}
\end{lem}
\begin{proof} See \cite{STL}.
\end{proof}
\begin{lem}\label{projnew}
Let $I=]a,b[$ and  $\Pi^{\alpha,\beta}_N$ be the $L^2_{\omega^{\alpha,\beta}}$-orthogonal projection. Then, for any $u \in H_{{\omega ^{\alpha ,\beta }},*}^m \left( I \right):= \left\{ {\left. u \right|\partial _t^lu \in L_{{\omega^{\alpha+l,\beta+l} }}^2\left( I  \right),0 \le l \le m} \right\}$, and $0\leq l\leq m\leq N+1$, there exists constant $C$  such that
\begin{equation}\label{57ttt}
{\left\| {\partial _t^l\left( {u - \Pi_N^{\alpha ,\beta }u} \right)} \right\|_{{\omega ^{\alpha + l,\beta + l}}}} \le C(\dfrac{(N-m+1)!}{(N-l+1)!})^{\frac{1}{2}}(N+m)^{\frac{l-m}{2}}{\left\| {\partial _t^mu} \right\|_{{\omega ^{\alpha + m,\beta + m}}}}.
\end{equation}
\end{lem}
\begin{proof} See \cite{STL}.
\end{proof}
For any $u$, $v$ on $I=[a, b]$, we set
\begin{equation*}
{\left( {u,v} \right)_{N,{\omega }}} = \sum\limits_{i = 0}^N {u({x_i})v({x_i}){w_i}},
\end{equation*}
where $x_i$ $(0 \leq i \leq N )$ are the Gauss, Gauss-Radau or Gauss-Lobatto quadrature nodes and $w_i$, $(0 \leq i \leq N )$, are the Gauss, Gauss-Radau or Gauss-Lobatto quadrature weights. The Gauss quadrature formulas imply that
\begin{equation} \label{aa991}
{\left( {u,v} \right)_{N,{\omega }}} ={\left( {u,v} \right)_{{\omega }}},\;\;\; if \;\;\; u.v \in {\mathbb{P}}_{2N+\delta},
\end{equation}
where $\delta = 1, 0,-1$ for Gauss, Gauss-Radau, and Gauss-Lobatto quadrature, respectively.

\begin{defin} \textup{\cite{Z. Wu}} Let $0<\alpha<1$ and $\Omega\subset \mathbb{R}^n$ be bounded. Then, $f\in C^{\alpha,\frac{\alpha}{2}}(\overline{\Omega}\times[0,T])$ if there exists a positive constant $C$ such that
\[|f(x_1,t_1)-f(x_2,t_2)|\leq C\Big(|x_1-x_2|^2+|t_1-t_2|\Big)^{\frac{\alpha}{2}},~\forall x_1,x_2\in\overline{\Omega},~\forall t_1,t_2\in[0,T].\] 
Furthermore, for any nonnegative integer $k$
\[ C^{2k+\alpha,k+\frac{\alpha}{2}}(\overline{\Omega}\times[0,T]):=\{ f\in C^{\alpha,\frac{\alpha}{2}}(\overline{\Omega}\times[0,T]): \partial_x^{\beta} \partial_t^{i} f\in C^{\alpha,\frac{\alpha}{2}}(\overline{\Omega}\times[0,T]),~~ |\beta |+2i\leq 2k\}.\]
\end{defin}

\noindent {\it\textbf{Proof of Lemma \ref{convpartial}}} Since $\dfrac{\partial^2 c}{\partial\rho^2}$  is a $C^1$-smooth function, therefore   for each $N^1\in\mathbb{N}$, there exists a polynomial $c_1^{N^1}$  such that
\[I_{N^1}^{0,0}c^{N^1}_1=c_1^{N^1},~\dfrac{\partial c_1^{N^1}}{\partial\rho }\left(0, t \right){\rm =0, }~c_1^{N^1}\left({\rm 1,} t \right){\rm =}0,~c_1^{N^1}\left(\rho ,0\right){\rm =}0, ~~0\leq t\leq T,\]
and
\begin{equation}\label{ABA1}
\lim_{N^1\rightarrow\infty} \Big(\|\dfrac{\partial c_1^{N^1}}{\partial\rho}-\dfrac{\partial c}{\partial\rho}\|^2_{\omega ^{0,0}}+\|I_{N^1}^{0,0}c-c_1^{N^1}\|^2_{\omega^{0,0}}+\|I_{N^1}^{0,0}\dfrac{\partial(c-c_1^{N^1})}{\rho\partial\rho}\|^2_{\omega^{0,0}}+\|I_{N^1}^{0,0}\dfrac{\partial^2(c-c_1^{N^1})}{\partial\rho^2}\|^2_{\omega^{0,0}}\Big)=0.
\end{equation}
Moreover, if $\dfrac{\partial^2 c}{\partial\rho^2}$ is a $C^m$-smooth function with respect to $\rho$ then it is easy to show that there exist $c_1^{N^1}$ and a positive constant $Q^*$ such that
\begin{equation}\label{ABA1kk}
\lim_{N^1\rightarrow\infty} \Big(\|\dfrac{\partial c_1^{N^1}}{\partial\rho}-\dfrac{\partial c}{\partial\rho}\|^2_{\omega ^{0,0}}+\|I_{N^1}^{0,0}c-c_1^{N^1}\|^2_{\omega^{0,0}}+\|I_{N^1}^{0,0}\dfrac{\partial(c-c_1^{N^1})}{\rho\partial\rho}\|^2_{\omega^{0,0}}+\|I_{N^1}^{0,0}\dfrac{\partial^2(c-c_1^{N^1})}{\partial\rho^2}\|^2_{\omega^{0,0}}\Big)\leq \dfrac{Q^*}{(N^1)^{2m}}.
\end{equation}
From (\ref{collo2}), we can conclude
\[( \Pi_{N^1}^{0,0}Lc^{N^1}_{n+1}-\Pi_{N^1}^{0,0}Lc^{N^1}_{n+1,1},\dfrac{\partial^2 (c^{N^1}_{n+1}- c^{N^1}_{n+1,1})}{\partial\rho^2})_{\omega^{0,0} } =\]
\[(I_{N^1}^{0,0}(\smallunderbrace{g^*_n-Lc^{N^1}_{n+1,1}}_{g_n^1}),\dfrac{\partial^2 (c^{N^1}_{n+1}- c^{N^1}_{n+1,1})}{\partial\rho^2})_{\omega ^{0,0}},\]
where $c^{N^1}_{n+1,1}(\rho):=c^{N^1}_1(\rho, t_{n+1})$ and for each $\phi\in C^2[0,1]$,  $L$ is defined as follows  \[L\phi :=\phi-\dfrac{a_0'D_{1}}{(R_{n+1}^{ap})^2}\dfrac{1}{{\rho }^{2}}\dfrac{\partial }{\partial\rho }\Big({\rho }^{2}\dfrac{\partial \phi}{\partial\rho }\Big)-t^*\dfrac{2(2{v}^{ap}(1,t_{n})-{v}^{ap}(1,t_{n-1}))\rho}{3R_{n+1}^{ap}}\dfrac{\partial \phi}{\partial\rho}.\]
 Therefore, we have
\[(c^{N^1}_{n+1}-c^{N^1}_{n+1,1},\dfrac{\partial^2(c^{N^1}_{n+1}- c^{N^1}_{n+1,1})}{\partial\rho^2})_{\omega ^{0,0}}=(I_{N^1}^{0,0}g_n^1,\dfrac{\partial^2 (c^{N^1}_{n+1}- c^{N^1}_{n+1,1})}{\partial\rho^2})_{\omega ^{0,0}}+\]
\[(\dfrac{a_0'D_1}{(R^{ap}_{n+1})^2}\dfrac{1}{{\rho^2 }}\dfrac{\partial }{\partial\rho }\Big({\rho^2 }\dfrac{\partial (c^{N^1}_{n+1}-c^{N^1}_{n+1,1})}{\partial\rho }\Big),\dfrac{\partial^2 (c^{N^1}_{n+1}- c^{N^1}_{n+1,1})}{\partial\rho^2})_{\omega ^{0,0}}+\]
\[(t^*\dfrac{2(2{v}^{ap}(1,t_{n})-{v}^{ap}(1,t_{n-1}))\rho}{3R_{n+1}^{ap}}\dfrac{\partial (c^{N^1}_{n+1}- c^{N^1}_{n+1,1})}{\partial\rho},\dfrac{\partial^2 (c^{N^1}_{n+1}- c^{N^1}_{n+1,1})}{\partial\rho^2})_{\omega ^{0,0}}=\]
\[\dfrac{a_0'D_1}{(R^{ap}_{n+1})^2}\|\dfrac{\partial^2 (c^{N^1}_{n+1}- c^{N^1}_{n+1,1})}{\partial \rho^2}\|^2_{\omega ^{0,0}}+\dfrac{a_0'D_1}{(R_{n+1}^{ap})^2}\|\dfrac{\partial (c^{N^1}_{n+1}- c_{n+1,1}^{N^1})}{\rho\partial \rho}\|^2_{\omega ^{0,0}}+\]
\[\dfrac{a_0'D_1}{(R^{ap}_{n+1})^2}(\dfrac{\partial (c^{N^1}_{n+1}- c^{N^1}_{n+1,1})}{\partial \rho})^2|_{\rho=1}+\]
\[(t^*\dfrac{2(2{v}^{ap}(1,t_{n})-{v}^{ap}(1,t_{n-1}))\rho}{3R_{n+1}^{ap}}\dfrac{\partial (c^{N^1}_{n+1}- c^{N^1}_{n+1,1})}{\partial\rho},\dfrac{\partial^2 (c^{N^1}_{n+1}- c^{N^1}_{n+1,1})}{\partial\rho^2})_{\omega ^{0,0}}+\]
\begin{equation}\label{khojad44}
(I_{N^1}^{0,0}g_n^1,\dfrac{\partial^2 (c^{N^1}_{n+1}- c^{N^1}_{n+1,1})}{\partial\rho^2})_{\omega^{0,0}}.
\end{equation}
Thus from \eqref{khojad44},  one can deduce that  
  \[\|\dfrac{\partial(c^{N^1}_{n+1}- c^{N^1}_{n+1,1})}{\partial\rho}\|^2_{\omega ^{0,0}}+\dfrac{a_0'D_1}{2(R^{ap}_{n+1})^2}\|\dfrac{\partial^2 (c^{N^1}_{n+1}- c^{N^1}_{n+1,1})}{\partial \rho^2}\|^2_{\omega ^{0,0}}+\dfrac{a_0'D_1}{2(R^{ap}_{n+1})^2}\|\dfrac{\partial (c^{N^1}_{n+1}- c^{N^1}_{n+1,1})}{\rho\partial \rho}\|^2_{\omega ^{0,0}}+\]
\[\dfrac{a_0'D_1}{(R^{ap}_{n+1})^2}(\dfrac{\partial (c^{N^1}_{n+1}- c^{N^1}_{n+1,1})}{\partial \rho})^2|_{\rho=1}\leq\]
\begin{equation}\label{imp156}
\Big|(I_{N^1}^{0,0}g_n^1,\dfrac{\partial^2 (c^{N^1}_{n+1}- c^{N^1}_{n+1,1})}{\partial\rho^2})_{\omega ^{0,0}}\Big|.
\end{equation}
In addition, from \eqref{ap00}, we have
 \begin{eqnarray}\label{imp276}
&\Big|(I_{N^1}^{0,0}g_n^1,\dfrac{\partial^2 (c^{N^1}_{n+1}- c^{N^1}_{n+1,1})}{\partial\rho^2})_{\omega^{0,0}}\Big|\leq\Big|(I_{N^1}^{0,0}(g^*_n-Lc_{n+1,1}^{N^1}),\dfrac{\partial^2 (c^{N^1}_{n+1}- c^{N^1}_{n+1,1})}{\partial\rho^2})_{\omega^{0,0}}\Big|\leq &\nonumber\\
&\Big|(c^{ap}_n-c^{N^1}_{n,1}-\dfrac{c^{ap}_{n-1}-c^{N^1}_{n-1,1}-c^{ap}_n+c^{N^1}_{n,1}}{3},\dfrac{\partial^2 (c^{N^1}_{n+1}- c^{N^1}_{n+1,1})}{\partial\rho^2})_{\omega^{0,0}}\Big|+&\nonumber\\
&\Big|\displaystyle(\sum_{k=0}^{n-1} (a'_k-a'_{k+1})\Big(\dfrac{D_{1}}{(R^{ap}_{n-k})^2}\dfrac{1}{{\rho }^{2}}\dfrac{\partial }{\partial\rho }({\rho }^{2}\dfrac{\partial c^{ap}_{n-k}-c^{N^1}_{n-k,1}}{\partial\rho })\Big),\dfrac{\partial^2 (c^{N^1}_{n+1}- c^{N^1}_{n+1,1})}{\partial\rho^2})_{\omega^{0,0}}\Big|+&\nonumber\\
&\Big|(\smallunderbrace{I_{N^1}^{0,0}L^2_nc_1^{N^1}-I_{N^1}^{0,0}f_n^*}_{g_n^2},\dfrac{\partial^2 (c^{N^1}_{n+1}- c^{N^1}_{n+1,1})}{\partial\rho^2})_{\omega^{0,0}}\Big|,&
\end{eqnarray}
where
\[L_n^2\phi=\phi_{n+1}-\phi_n+\dfrac{\phi_{n-1}-\phi_n}{3}-\]
\[t^*\dfrac{2(2{v}^{ap}(1,t_{n})-{v}^{ap}(1,t_{n-1}))\rho}{3R_{n+1}^{ap}}\dfrac{\partial \phi_{n+1}}{\partial\rho}-\dfrac{a_0'D_{1}}{(R^{ap}_{n+1})^2}\dfrac{1}{{\rho }^{2}}\dfrac{\partial }{\partial\rho }\left({\rho }^{2}\dfrac{\partial \phi_{n+1}}{\partial\rho }\right)+\]
\[\displaystyle\sum_{k=0}^{n-1} (a'_k-a'_{k+1})\Big(\dfrac{D_{1}}{(R^{ap}_{n-k})^2}\dfrac{1}{{\rho }^{2}}\dfrac{\partial }{\partial\rho }\left({\rho }^{2}\dfrac{\partial \phi_{n-k}}{\partial\rho }\right)\Big).\]
Therefore, using \eqref{imp276}, Cauchy--Schwarz inequality  and \eqref{ak0llll}, for each positive $\epsilon$ and $\epsilon_1$, we can deduce that  
 \begin{eqnarray}\label{imp256}
&\Big|(I_{N^1}^{0,0}g_n^1,\dfrac{\partial^2 (c^{N^1}_{n+1}- c^{N^1}_{n+1,1})}{\partial\rho^2})_{\omega^{0,0}}\Big|\leq \Big|({g_n^2},\dfrac{\partial^2 (c^{N^1}_{n+1}- c^{N^1}_{n+1,1})}{\partial\rho^2})_{\omega^{0,0}}\Big|+ &\nonumber\\
&\displaystyle\sum_{k=0}^{n-1}\dfrac{(a'_k-a'_{k+1})D_1}{\epsilon(R^{ap}_{n-k})^2}\| \dfrac{\partial^2(c^{ap}_{n-k}-c_{n-k,1}^{N^1})}{\partial \rho^2}\|_{\omega^{0,0}}^2+{\epsilon}\displaystyle\sum_{k=0}^{n-1}\dfrac{(a'_k-a'_{k+1})D_1}{(R^{ap}_{n-k})^2}\|\dfrac{\partial^2 (c^{N^1}_{n+1}- c^{N^1}_{n+1,1})}{\partial\rho^2}\|^2_{\omega^{0,0}}+&\nonumber\\
&\dfrac{1}{2\epsilon_1}\|\dfrac{\partial (c^{ap}_n-c^{N^1}_{n,1})}{\partial\rho}\|^2_{\omega^{0,0}}+\displaystyle\sum_{k=0}^{n-1}\dfrac{(a'_k-a'_{k+1})D_1}{\epsilon(R^{ap}_{n-k})^2}\| \dfrac{2}{\rho}\dfrac{\partial(c^{ap}_{n-k}-c_{n-k,1}^{N^1})}{\partial \rho}\|_{\omega^{0,0}}^2+&\nonumber\\
&\dfrac{\epsilon_1}{2}\|\dfrac{\partial (c^{N^1}_{n+1}- c^{N^1}_{n+1,1})}{\partial\rho}\|^2_{\omega^{0,0}}+\dfrac{1}{3}\Big|(\dfrac{\partial (c^{ap}_{n-1}-c^{N^1}_{n-1,1}-c^{ap}_n+c^{N^1}_{n,1})}{\partial\rho},\dfrac{\partial (c^{N^1}_{n+1}- c^{N^1}_{n+1,1})}{\partial\rho})_{\omega^{0,0}}\Big|.&
\end{eqnarray}
Therefore, from  \eqref{imp156} and \eqref{imp256}, we deduce that  there exists a positive constant $K_2$ such that
  \[\dfrac{1}{2}\|\dfrac{\partial(c^{N^1}_{n+1}- c^{N^1}_{n+1,1})}{\partial\rho}\|^2_{\omega ^{0,0}}+\dfrac{D_1a_0'}{3(R^{ap}_{n+1})^2}\|\dfrac{\partial^2 (c^{N^1}_{n+1}- c^{N^1}_{n+1,1})}{\partial \rho^2}\|^2_{\omega ^{0,0}}+\dfrac{D_1a_0'}{2(R^{ap}_{n+1})^2}\|\dfrac{\partial (c^{N^1}_{n+1}- c^{N^1}_{n+1,1})}{\rho\partial \rho}\|^2_{\omega ^{0,0}}+\]
\[\dfrac{D_1a_0'}{(R^{ap}_{n+1})^2}(\dfrac{\partial (c^{N^1}_{n+1}- c^{N^1}_{n+1,1})}{\partial \rho})^2|_{\rho=1}\leq\]
\[ \dfrac{1}{2}\|\dfrac{\partial (c^{ap}_n-c^{N^1}_{n,1})}{\partial\rho}\|^2_{\omega^{0,0}}+{K_2}\Big(\Big|({g_n^2},\dfrac{\partial^2 (c^{N^1}_{n+1}- c^{N^1}_{n+1,1})}{\partial\rho^2})_{\omega^{0,0}}\Big|+\displaystyle\sum_{k=0}^{n-1}\dfrac{(a'_k-a'_{k+1})D_1}{(R^{ap}_{n-k})^2}\| \dfrac{2}{\rho}\dfrac{\partial(c^{ap}_{n-k}-c_{n-k,1}^{N^1})}{\partial \rho}\|_{\omega^{0,0}}^2+\]
\[\displaystyle\sum_{k=0}^{n-1}\dfrac{(a'_k-a'_{k+1})D_1}{(R^{ap}_{n-k})^2}\| \dfrac{\partial^2(c^{ap}_{n-k}-c_{n-k,1}^{N^1})}{\partial \rho^2}\|_{\omega^{0,0}}^2\Big)+\dfrac{1}{3}\Big|(\dfrac{\partial (c^{ap}_{n-1}-c^{N^1}_{n-1,1}-c^{ap}_n+c^{N^1}_{n,1})}{\partial\rho},\dfrac{\partial (c^{N^1}_{n+1}- c^{N^1}_{n+1,1})}{\partial\rho})_{\omega^{0,0}}\Big|.\]
 Hence, from \textbf{A}, \eqref{A1} and \eqref{A1**}, one can conclude that there exists positive $K_3$ such that
  \[\dfrac{1}{2}\|\dfrac{\partial(c^{N^1}_{n+1}- c^{N^1}_{n+1,1})}{\partial\rho}\|^2_{\omega ^{0,0}}+\dfrac{D_1a_0'}{4(R^{ap}_{n+1})^2}\|\dfrac{\partial^2 (c^{N^1}_{n+1}- c^{N^1}_{n+1,1})}{\partial \rho^2}\|^2_{\omega ^{0,0}}+\dfrac{D_1a_0'}{2(R^{ap}_{n+1})^2}\|\dfrac{\partial (c^{N^1}_{n+1}- c^{N^1}_{n+1,1})}{\rho\partial \rho}\|^2_{\omega ^{0,0}}+\]
\[\dfrac{D_1a_0'}{(R^{ap}_{n+1})^2}(\dfrac{\partial (c^{N^1}_{n+1}- c^{N^1}_{n+1,1})}{\partial \rho})^2|_{\rho=1}\leq\]
\[\dfrac{1}{2}\|\dfrac{\partial (c^{ap}_n-c^{N^1}_{n,1})}{\partial\rho}\|^2_{\omega^{0,0}}+K_3\Big( \displaystyle\sum_{k=0}^{n-1}\dfrac{(a'_k-a'_{k+1})D_1}{(R^{ap}_{n-k})^2}\| \dfrac{2}{\rho}\dfrac{\partial(c^{ap}_{n-k}-c_{n-k,1}^{N^1})}{\partial \rho}\|_{\omega^{0,0}}^2+\]
\[\displaystyle\sum_{k=0}^{n-1}\dfrac{(a'_k-a'_{k+1})D_1}{(R^{ap}_{n-k})^2}\| \dfrac{\partial^2(c^{ap}_{n-k}-c_{n-k,1}^{N^1})}{\partial \rho^2}\|_{\omega^{0,0}}^2+\]
\[\dfrac{1}{(t^*)^{1-\alpha}}\|2t^*f\left(c_n,p_n,q_n\right)-t^*f\left(c_{n-1},p_{n-1},q_{n-1}\right)+E_{ t}^c-2t^*f\left(c^{ap}_{n},p^{ap}_{n},q^{ap}_{n}\right)+t^*f\left(c^{ap}_{n-1},p^{ap}_{n-1},q^{ap}_{n-1}\right)\|^2_{\omega^{0,0}}+\]
\begin{equation}\label{tyu} 
\Big|(L_n^2c_{1}^{N^1}-L_n^*c_,\dfrac{\partial^2 (c^{N^1}_{n+1}- c^{N^1}_{n+1,1})}{\partial\rho^2})_{\omega^{0,0}}\Big|\Big)+\dfrac{1}{3}\Big|(\dfrac{\partial (c^{ap}_{n-1}-c^{N^1}_{n-1,1}-c^{ap}_n+c^{N^1}_{n,1})}{\partial\rho},\dfrac{\partial(c^{N^1}_{n+1}- c^{N^1}_{n+1,1})}{\partial\rho})_{\omega ^{0,0}}\Big|,
\end{equation}
where
\[L_n^*c=c_{n+1}-c_n-\dfrac{c_{n}-c_{n-1}}{3}-t^*\dfrac{2(2{v}(1,t_{n})-{v}(1,t_{n-1}))\rho}{3R_{n+1}}\dfrac{\partial c_{n+1}}{\partial\rho}-\]
\[\dfrac{a_0'D_{1}}{(R_{n+1})^2}\dfrac{1}{{\rho }^{2}}\dfrac{\partial }{\partial\rho }\left({\rho }^{2}\dfrac{\partial c_{n+1}}{\partial\rho }\right)+\displaystyle\sum_{k=0}^{n-1} (a'_k-a'_{k+1})\Big(\dfrac{D_{1}}{(R_{n-k})^2}\dfrac{1}{{\rho }^{2}}\dfrac{\partial }{\partial\rho }\left({\rho }^{2}\dfrac{\partial c_{n-k}}{\partial\rho }\right)\Big).\]
Therefore, from Lemma \ref{gronwall} (discrete Gronwall lemma), \eqref{TM6KK}, \eqref{TM6KK**}, \eqref{A1}, \eqref{A1**}, \eqref{ABA1} and  \eqref{tyu}, we can conclude that there exists positive $K_4$ such that
 \[\|\dfrac{\partial(c^{N^1}_{n+1}- c_{n+1,1}^{N^1})}{\partial\rho}\|^2_{\omega ^{0,0}}\leq \|\dfrac{\partial(c^{ap}_{n}- c_{n,1}^{N^1})}{\partial\rho}\|^2_{\omega ^{0,0}}+\]
\[ {K_4}\Big(({t^*})^{1+\alpha}\sum_{k=0}^n(\dfrac{1}{3})^{\frac{k}{2}}\Big(\|{c_{k,1}^{N^1}}-{c^{ap}_k }\|^2_{\omega^{0,0}}+\|{{p_k }}-{{p^{ap}_k }}\|^2_{\omega^{0,0}}+\|q_k -q^{ap}_k\|^2_{\omega^{0,0}}+\|d_k -d^{ap}_k\|^2_{\omega^{0,0}}+\|R_k -R^{ap}_k\|^2_{\omega^{0,0}}\Big)+\]
\[\|e_t^c\|^2_{\omega^{0,0}}+(t^*)K^*(N^1)\Big),\]
where  
\[\lim_{N^1\rightarrow\infty}K^*(N^1)=0,~~\|e_t^c\|_{\infty}\leq (t^*)^{2+\frac{1-\alpha}{2}}.\]
In addition, if $\dfrac{\partial^2 c}{\partial\rho^2}$ is a $C^m$-smooth function with respect to $\rho$ then from \eqref{ABA1kk}, we can show 
\[\|\dfrac{\partial(c^{N^1}_{n+1}- c_{n+1,1}^{N^1})}{\partial\rho}\|^2_{\omega ^{0,0}}\leq \|\dfrac{\partial(c^{ap}_{n}- c_{n,1}^{N^1})}{\partial\rho}\|^2_{\omega ^{0,0}}+ \]
\[K_4\Big(({t^*})^{1+\alpha}\sum_{k=0}^n(\dfrac{1}{3})^{\frac{k}{2}}\Big(\|{c_{k,1}^{N^1}}-{c^{ap}_k }\|^2_{\omega^{0,0}}+\|{{p_k }}-{{p^{ap}_k }}\|^2_{\omega^{0,0}}+\|q_k -q^{ap}_k\|^2_{\omega^{0,0}}+\|d_k -d^{ap}_k\|^2_{\omega^{0,0}}+\|R_k -R^{ap}_k\|^2_{\omega^{0,0}}\Big)+\]
\[\|e_ t^c\|^2_{\omega^{0,0}}+\dfrac{t^*}{({N^1})^{2m}}\Big).\]
\qed

\noindent{\it\textbf{Proof of Theorem \ref{mainthm}}} 
 Employing    \eqref{YABss}--\eqref{TM6KK}, \eqref{imp444}, and \eqref{ap33}--\eqref{TM6KK**}, one can conclude that there exist positive $M$ and $M^*$ such that 
 \begin{equation}\label{AL23345}
\mathfrak{E}_{k+1}\leq\mathfrak{E}_{k-1}+M t^* (\mathfrak{E}_k+\mathfrak{e}_k)+\|E_ t\|_{\infty}, ~~\|E_ t\|_{\infty}\leq M^*( t^*)^3,
\end{equation}
where 
\begin{equation}\label{defpqd}
\mathfrak{E}_k=\| p^{ap}_k-p_k\|_{\infty}+\|q^{ap}_k-q_k\|_{\infty}+\|d^{ap}_k-d_k\|_{\infty}+|R^{ap}_{k}-R_{k}|,
\end{equation}
and 
\begin{equation}\label{defcw}
\mathfrak{e}_k=\| \dfrac{\partial(c^{ap}_k-c_k)}{\partial\rho}\|_{\omega^{0,0}}+\|\dfrac{\partial(w^{ap}_k-w_k)}{\partial\rho}\|_{\omega^{0,0}}.
\end{equation}
Using  \eqref{AL23345}--\eqref{defcw},  one can conclude that  
\[\mathfrak{E}_{n+1}\leq\mathfrak{E}_{n-1}+M t^* (\mathfrak{E}_n+\mathfrak{e}_n)+\|E_ t\|_{\infty}\leq\]
\[(1+M t^*)\max\{\mathfrak{E}_{n-1}, \mathfrak{E}_{n}\}+M t^*\mathfrak{e}_n+\|E_ t\|_{\infty}\leq\] 
\begin{equation}\label{imp777}
\dfrac{(1+M t^*)^m-1}{M t^*}(\|E_ t\|_{\infty}+M t^*\max_{k=0,1,\cdots,n}\{\mathfrak{e}_{k}\}).
\end{equation}
 On the other hand, employing \eqref{HM1}--\eqref{HM2}, we get
\[\max_{k=0,1,\cdots,n}\{(\mathfrak{e}_{k})^2\}\leq\] \begin{equation}\label{imp888}
 \dfrac{M_1}{t^*}\Big((t^*)^{1+\alpha}\max_{k=0,1,\cdots,n}\{(\mathfrak{E}_{k})^2\}+\|E_t^1\|^2_{\infty}\Big)+M_1K_3^*(N^1),
\end{equation}
where $M_1$ is a positive constant and
\[\|E_t^1\|_{\infty}\leq (t^*)^{2+\frac{1-\alpha}{2}},~~K_3^*(N^1)=K_1^*(N^1)+K_2^*(N^1).\]
From \eqref{imp777} and \eqref{imp888}, we deduce that there exists positive $C_4^*$ such that
\begin{equation}\label{imp1010}
\max_{k=0,1,\cdots,n+1}\{\mathfrak{E}_{k}\}\leq C_4^* e^{MT}\Big(( t^*)^{2-\frac{\alpha}{2}}+  (K^*_3(N^1))^{\frac{1}{2}}\Big).
\end{equation}
Therefore, employing \eqref{imp888} and \eqref{imp1010}, we conclude that there exists positive constant $C_5^*$ such that
 \begin{equation}\label{imp1111}
\max_{k=0,1,\cdots,n+1}\{\mathfrak{e}_{k}\}\leq C_5^* e^{MT}\Big(( t^*)^{2-\frac{\alpha}{2}}+   (K^*_3(N^1))^{\frac{1}{2}}\Big).
\end{equation}
Finally, employing \eqref{HM1} and \eqref{HM1**KF}, we can get the desired results.
\qed

 {\it\noindent\textbf{Proof of Theorem \ref{stability}}}
If we solve the perturbed problem \eqref{YAS88900}-\eqref{HM34588900}  using the presented method, employing    \eqref{YABss}--\eqref{TM6KK}, \eqref{imp444}, and \eqref{ap33}--\eqref{TM6KK**} for perturbed problem, one can conclude that there exist positive $M_2$ and $M_2^*$ such that 
 \begin{equation}\label{AL233453}
\mathfrak{E}^p_{k+1}\leq\mathfrak{E}^p_{k-1}+M_2 t^* (\mathfrak{E}^{p}_k+\mathfrak{e}^p_k)+ t^*\epsilon_1+\|E^p_ t\|_{\infty}, ~~\|E^p_ t\|_{\infty}<M_2^*( t^*)^3,
\end{equation}
where 
\begin{equation}\label{defpqd3}
\mathfrak{E}^p_k=\| p^{ap,p}_k-p_k\|_{\infty}+\|q^{ap,p}_k-q_k\|_{\infty}+\|d^{ap,p}_k-d_k\|_{\infty}+|R^{ap,p}_{k}-R_{k}|,
\end{equation}
and 
\begin{equation}\label{defcw3}
\mathfrak{e}^p_k=\| \dfrac{\partial(c^{ap,p}_k-c_k)}{\partial\rho}\|_{\omega^{0,0}}+\|\dfrac{\partial(w^{ap,p}_k-w_k)}{\partial\rho}\|_{\omega^{0,0}}.
\end{equation}
Using  \eqref{AL233453}--\eqref{defcw3},  one can conclude that  
\[\mathfrak{E}^p_{n+1}\leq\mathfrak{E}^p_{n-1}+M_2 t^* (\mathfrak{E}^p_n+\mathfrak{e}^p_n)+ t^*\epsilon_1+\|E^p_ t\|_{\infty}\leq\]
\[(1+M_2 t^*)\max\{\mathfrak{E}^p_{n-1}, \mathfrak{E}^p_{n}\}+M_2 t^*\mathfrak{e}^p_n+ t^*\epsilon_1+\|E^p_ t\|_{\infty}\leq\] 
\begin{equation}\label{imp7773}
\dfrac{(1+M_2 t^*)^m-1}{M_2 t^*}( t^*\epsilon_1+\|E^p_ t\|_{\infty}+M_2 t^*\max_{k=0,1,\cdots,n}\{\mathfrak{e}^p_{k}\}).
\end{equation}
Also, from \eqref{HM1}--\eqref{HM2}, we conclude that
 \begin{equation}\label{imp8883}
\max_{k=0,1,\cdots,n}\{\mathfrak{e}^{p}_{k}\}\leq M_3\Big((t^*)^{\frac{\alpha}{2}}(\max_{k=0,1,\cdots,n}\{\mathfrak{E}^p_{k}\})+\epsilon_1+K_5^*(N^1)+\|E_t^1\|_{\infty}\Big),
\end{equation}
where $M_3$ is a positive constant and
\[\|E_t^1\|_{\infty}\leq C_6(t^*)^{2-\frac{\alpha}{2}},~~K_5^*(N^1)=(K_1^*(N^1))^{\frac{1}{2}}+(K_2^*(N^1))^{\frac{1}{2}}.\] 
From \eqref{imp7773} and \eqref{imp8883}, we deduce that there exists positive $C_6^*$ such that
\begin{equation}\label{imp10103}
\max_{k=0,1,\cdots,n+1}\{\mathfrak{E}^p_{k}\}\leq C_6^* e^{M_2T}\Big(( t^*)^{2-\frac{\alpha}{2}}+  K^*_5(N^1)+\epsilon_1\Big).
\end{equation}
Therefore, employing \eqref{imp8883} and \eqref{imp10103}, we conclude that there exists positive constant $C_7^*$ such that
 \begin{equation}\label{imp11113}
\max_{k=0,1,\cdots,n+1}\{\mathfrak{e}^p_{k}\}\leq C_7^* e^{M_2T}\Big(( t^*)^{2-\frac{\alpha}{2}}+ K^*_5(N^1)+\epsilon_1\Big).
\end{equation}
\qed


\begin{thebibliography}{amsplain}
\bibitem {glio1} Y. Kim,   Regulation of cell proliferation and
migration in glioblastoma: new therapeutic
approach, Front. Oncol. 3 (2013) 359--371.

\bibitem {glio} J. C. L. Alfonso, K. Talkenberger, M. Seifert, B. Klink, A. Hawkins-Daarud, K. R. Swanson, H. Hatzikirou, A. Deutsch,  The biology and mathematical modelling
of glioma invasion: a review, J. R. Soc. Interface (2017)  DOI: 10.1098/rsif.2017.0490. 

\bibitem {phyll} L. Duman,  N.S. Gezer,  P. Balci, c. Altay,  I. Basara,  M.G. Durak, A.I.
 Sevinc, Differentiation
between phyllodes tumors and fibroadenomas based on mammographic
sonographic and MRI features, Breast Care  11 (2016) 123--127.


\bibitem {vas1} T. L. Jackson, Vascular tumor growth and treatment: Consequences of polyclonality, competition and dynamic vascular support, J. Math. Biol. 44 (2002) 201--226. 
\bibitem {J.H.} J.  Zhao,  A parabolic-hyperbolic free boundary problem modeling tumor growth with drug application, Electron. J. Differ. Eq.  2010 (2010) 1--18. 

\bibitem {vas2}T. L. Jackson, H. M. Byrne, A mathematical model to study the effects of drug resistance and vasculature on the response of solid tumors to chemotherapy, Math. Biosci. 164 (2000), 17--38.

\bibitem {M. Chen}Y. Tao, M. Chen, An elliptic-hyperbolic free boundary problem modelling cancer therapy, Nonlinearity, 19 (2006), 419--440. 



\bibitem {khai} D. Khaitan, S. Chandna, M.B. Arya, B.S. Dwarakanath,  Establishment and characterization of
multicellular spheroids from a human glioma cell line: implications for tumor therapy, J. Transl. Med.
4 (2006) 12--25 

\bibitem {Y. Tao}Y. Tao,  A free boundary problem modeling the cell cycle and cell movement in multicellular tumor spheroids, J. Diff. Eq.  247 (2009) 49--68.

\bibitem{ther1}  C. Geng, H. Paganetti,  C. Grassberger,  Prediction of treatment response for combined chemo and radiation therapy for non-small cell lung cancer patients using a bio-mathematical model, Sci Rep. 7 (2017) 13542 doi:10.1038/s41598-017-13646-z.

\bibitem{ther2} K. Abernathy, Z. Abernathy, K. Brown, C. Burgess, R. Hoehne, Global dynamics of a colorectal cancer treatment model with cancer stem cells, Heliyon.  3 (2017)e00247 doi: 10.1016/j.heliyon.2017.e00247. 
\bibitem{ther3} H. Enderling, M.A. Chaplain,  Mathematical
modeling of tumor growth and treatment. Curr.
Pharm. Des. 20 (2014) 4934--4940.



\bibitem{fansari}  F. Ansarizadeh, M. Singh,  D. Richards, Modelling of tumor cells
regression in response to chemotherapeutic treatment, Appl. Math. Model.
48  (2017) 96-112




 
\bibitem{Ates}  I. Ates and P. A. Zegeling, A homotopy perturbation method for fractional
order advection-diffusion-reaction boundary-value problems, Appl. Math.
Model. 47 (2017) 425--441 


\bibitem{Sohail} A. Sohail, S. Arshad, S. Javed, K. Maqbool, Numerical analysis of fractional-order tumor model, Int. J. Biomath. 8 (2015) 1550069


\bibitem{Veeresha}   P. Veeresha, D.G. Prakasha,  H.M.  Baskonus, New numerical surfaces to the mathematical model of cancer
chemotherapy effect in Caputo fractional derivatives, Chaos 29 (2019)  013119. doi:10.1063/1.5074099.


\bibitem{Qu} A. Quarteroni, A. Valli, Numerical Approximation of Partial Differential Equations, Springer,
Berlin, 1997.








 


\bibitem {CHQZ} C. Canuto, M.Y. Hussaini, A. Quarteroni, T. A. Zang,   Spectral Methods Fundamentals in
Single Domains,  Springer, Berlin, 2006.
\bibitem {STL} J. Shen, T. Tang, L. Wang, Spectral Methods, Algorithms, Analysis and Applications,  Springer-Verlag Berlin Heidelberg, 2011.

\bibitem{Z. Wu}  Z. Wu, J. Yin,  C. Wang,   {Elliptic and Parabolic Equations,} World Scientific, Singapore, 2006. 

\bibitem{impfrac} Y.M. Lin, C.J. Xu, Finite difference/spectral approximations for the time-fractional diffusion
equation, J. Comput. Phys. 225 (2007) 1533-1552


\bibitem{Huang}  C. Huang,  Z. Zhang, Q. Song,  Spectral methods for substantial fractional differential
equations, J. Sci. Comput. 74 (2018) 1554--1574


\end{thebibliography}
\end{document}